\documentclass[12pt]{article}

\usepackage{subcaption}

\usepackage{float}

\usepackage{multirow}

\usepackage{tikz}

\usetikzlibrary{calc,backgrounds,arrows,matrix}

\usepackage{enumerate}

\usepackage{color}
\definecolor{lightblue}{rgb}{0,0.2,0.5}

\usepackage{xcolor}
\definecolor{ForestGreen}{RGB}{34,139,34}
\definecolor{mauve}{rgb}{0.7,0,0.43}
\definecolor{dkgreen}{rgb}{0,0.6,0}
\definecolor{darkgreen}{rgb}{0,0.6,0}
\definecolor{darkorange}{rgb}{1.0, 0.55, 0.0}
\definecolor{lightblue}{rgb}{0,0.2,0.5}
\definecolor{blue1}{rgb}{0,0.1,0.9}

\usepackage[colorlinks=true, urlcolor=lightblue,linkcolor=lightblue, citecolor=lightblue]{hyperref}

\usepackage{listings}
\lstdefinelanguage{Maple}{
    morekeywords={proc, if, return, map, op, int, for, do, local, nops, convert, end},
    sensitive=false, 
    morecomment=[l]{//}, 
    morecomment=[s]{/*}{*/}, 
    morestring=[b]" 
} 

\usepackage{amssymb,amsmath}

\usepackage{graphicx}

\usepackage{graphics}

\usepackage{tikz}
\usetikzlibrary{trees}

\DeclareMathAlphabet{\eufrak}{U}{}{}{}
\SetMathAlphabet\eufrak{normal}{U}{euf}{m}{n}
\SetMathAlphabet\eufrak{bold}{U}{euf}{b}{n}

\allowdisplaybreaks

 \def\qu{{\mathord{\mathbb Z}}}

 \def\T{T} 

 \def\inte{{\mathord{\mathbb R}}}

 \def\inte{{\mathord{\mathbb N}}}

 \def\sZZ{{\rm Z\kern-.45em{}Z}}

 \def\sQQ{{\kern 0.27em \vrule height1.45ex width0.03em depth0em
           \kern-0.30em \rm Q}}
 \def\qu{{\mathchoice
         {\sQQ}
         {\sQQ}
   {\kern 0.225em \vrule height1.05ex width0.025em depth0em \kern-0.25em \rm Q}
   {\kern 0.180em \vrule height0.78ex width0.020em depth0em \kern-0.20em \rm Q}
         }}
 \def\sGG{{\kern 0.27em \vrule height1.45ex width0.03em depth0em
           \kern-0.30em \rm G}}
 \def\gg{{\mathchoice
         {\sGG}
         {\sGG}
   {\kern 0.225em \vrule height1.05ex width0.025em depth0em \kern-0.25em \rm G}
   {\kern 0.180em \vrule height0.78ex width0.020em depth0em \kern-0.20em \rm G}
         }}

 \newtheorem{prop}{Proposition}[section]
 \newtheorem{lemma}[prop]{Lemma}
 \newtheorem{definition}[prop]{Definition}
 
 \newtheorem{theorem}[prop]{Theorem}
 \newtheorem{remark}[prop]{Remark}

\numberwithin{equation}{section}

 \def\P{{\mathord{\mathbb P}}}

\newcommand{\re}{\mathrm{e}}

 \newcounter{hyp}
\newenvironment{Proof}{\removelastskip\par\medskip \noindent{\em Proof.} \rm}{\penalty-20\null\hfill$\square$\par\medbreak}

\def\bprf{\begin{Proof}}
\def\nprf{\end{Proof}}
\def\bdes{\begin{description}}
\def\ndes{\end{description}}

\newtheorem{thm}{Theorem}[section]

\renewcommand{\leq}{\leqslant}

\renewcommand{\geq}{\geqslant}

\def\bdef{\begin{defn}}
\def\ndef{\end{defn}}
\def\bthm{\begin{thm}}
\def\nthm{\end{thm}}
\def\bprop{\begin{prop}}
\def\nprop{\end{prop}}
\def\brmk{\begin{remark}}
\def\nrmk{\end{remark}}
\def\bexa{\begin{exa}}
\def\nexa{\end{exa}}
\def\blem{\begin{lem}}
\def\nlem{\end{lem}}
\def\bcor{\begin{cor}}
\def\ncor{\end{cor}}
\def\bexe{\begin{exe}}
\def\nexe{\end{exe}}

\newcommand{\E}{\mathbb{E}}

\newcommand{\real}{\mathbb{R}}

\def\og{\leavevmode\raise.3ex
     \hbox{$\scriptscriptstyle\langle\!\langle$~}}
\def\fg{\leavevmode\raise.3ex
     \hbox{~$\!\scriptscriptstyle\,\rangle\!\rangle$}~}

\title{\Huge
 A fully nonlinear Feynman-Kac formula with derivatives of arbitrary orders 
} 

\author{
 Jiang Yu Nguwi\footnote{\href{mailto:nguw0003@e.ntu.edu.sg}{nguw0003@e.ntu.edu.sg}
 }
 \qquad
  Guillaume Penent\footnote{\href{mailto:PENE0001@e.ntu.edu.sg}{pene0001@e.ntu.edu.sg}}
  \qquad
      Nicolas Privault\footnote{
\href{mailto:nprivault@ntu.edu.sg}{nprivault@ntu.edu.sg}
}
  \\
\small
Division of Mathematical Sciences
\\
\small
School of Physical and Mathematical Sciences
\\
\small
Nanyang Technological University
\\
\small
21 Nanyang Link, Singapore 637371
}

\allowdisplaybreaks

\usepackage{empheq}

\usepackage{bbm}

\makeatletter
\newcommand*\rel@kern[1]{\kern#1\dimexpr\macc@kerna}
\newcommand*\widebar[1]{
  \begingroup
  \def\mathaccent##1##2{
    \rel@kern{0.8}
    \overline{\rel@kern{-0.8}\macc@nucleus\rel@kern{0.2}}
    \rel@kern{-0.2}
  }
  \macc@depth\@ne
  \let\math@bgroup\@empty \let\math@egroup\macc@set@skewchar
  \mathsurround\z@ \frozen@everymath{\mathgroup\macc@group\relax}
  \macc@set@skewchar\relax
  \let\mathaccentV\macc@nested@a
  \macc@nested@a\relax111{#1}
  \endgroup
}
\makeatother

\begin{document}

\maketitle

\baselineskip0.6cm

\vspace{-0.6cm}

\begin{abstract}
 We present an algorithm for the numerical solution of
 nonlinear parabolic partial differential equations. 
 This algorithm extends the classical Feynman-Kac formula
 to fully nonlinear partial differential equations, by using random trees
 that carry information on nonlinearities on their branches. 
 It applies to functional, non-polynomial nonlinearities that
 are not treated by standard branching arguments, 
 and deals with derivative terms of arbitrary orders.
 A Monte Carlo numerical implementation is provided. 
\end{abstract}



\baselineskip0.7cm

\section{Introduction}

The objective of this paper is to provide probabilistic representations
 for the solutions of fully nonlinear parabolic partial
differential equations involving higher order derivatives, of the form
\begin{equation}
\label{eq:1}
\begin{cases}
  \displaystyle
  \partial_t u (t,x)
 + \frac{1}{2} \partial^2_x u (t,x)
  + 
  f\big(u(t,x) , \partial_x u (t,x), \ldots , \partial^n_x u (t,x) \big)
  = 0,
  \medskip
  \\
u(T,x) = \phi(x), \qquad (t,x) \in [0,T]\times \real,
\end{cases}
\end{equation}
 $n\geq 0$, where $\partial_x^2$ is the standard Laplacian on $\real$ and
 $f(x,y,z_1,\ldots , z_m)$ is a 
   smooth functional nonlinearity involving derivatives of arbitrary orders.
Probabilistic representations for the solutions of
first order nonlinear partial differential equations (PDEs) of the form 
\begin{equation} 
\nonumber 
\left\{ 
\begin{array}{l} 
  \displaystyle
  \partial_t u (t,x) 
 + \frac{1}{2} \partial^2_x u(t,x) 
 + f \big( t , x , u(t,x), \partial_x u (t,x) \big) = 0
\\ 
\\ 
 u(T,x) = \phi (x), 
 \qquad x \in \real^d,  
\end{array} 
\right. 
\end{equation} 
 can be obtained using backward stochastic differential 
 equations (BSDEs) \cite{peng2}, \cite{pardouxpeng}, by
 representing $u(t,x)$ as 
 $u(t,x) = Y_t^{t,x}$, $(t,x)\in [0,T]\times \real$,
 where $(Y_s^{t,x})_{t \leq s \leq T}$ is the 
 solution of the backward stochastic differential equation 
\begin{equation} 
\nonumber 
\left\{ 
\begin{array}{l} 
 dY_s^{t,x} = - f\big( s , X_s^{t,x} , Y_s^{t,x} , Z_s^{t,x} \big) ds 
 + 
 Z_s^{t,x} dX^{t,x}_s, 
 \qquad 
 0 \leq t \leq s \leq T,  
\\ 
\\ 
 Y_T^{t,x} = \phi ( X_T^{t,x} ) 
, 
\end{array} 
\right. 
\end{equation} 
 with random anticipating terminal condition, 
 where $(X^{t,x}_s)_{t \leq s \leq T}$ is a 
 standard Brownian motion started at $x\in \real$ at 
 time $t\in [0,T]$. 
 This method also extends to second order fully nonlinear PDEs of
 the form
$$ 
\left\{ 
\begin{array}{l} 
  \displaystyle
  \partial_t u (t,x) 
 +  
 f \big( t , x , u(t,x),  \partial_x u (t,x) ,
 \partial^2_x u (t,x) \big) = 0 
\medskip 
\\ 
 u(T,x) = \phi (x), 
 \qquad x \in \real^d,  
\end{array} 
\right. 
$$ 
 using second order backward stochastic differential
 equations \cite{touzi}, \cite{soner},
 see \cite{han2018solving},
 \cite{germain},
 \cite{lefebvre} 
 for deep learning implementations. 
 
\medskip 

On the other hand, stochastic diffusion branching mechanisms
for the representation of solutions of partial differential
equations have been introduced in \cite{skorohodbranching},
and extended to branching Markov processes in \cite{inw}.
Branching diffusions have also been applied to 
give a probabilistic representation of the solutions
of the Kolmogorov-Petrovskii-Piskunov (KPP) equation in \cite{hpmckean},
and to more general PDEs with polynomial nonlinearities in 
\cite{henry-labordere2012}, 
\cite{laborderespa}, see also \cite{chakraborty} for existence of solutions
of parabolic PDEs with 
power series nonlinearities. 

\medskip

\indent
In \cite{labordere}, this branching argument has been applied 
to polynomial gradient nonlinearities using branching trees.
In this approach, branches corresponding to gradient terms 
 are identified by marks and associated random weights
 which are used in Malliavin integration by parts,
 see also \cite{labordere2}
 for an application to semilinear and higher-order hyperbolic PDEs. 
 In \cite{fahim}, see also \cite{tanxiaolu}, \cite{guowenjie}, \cite{kong}, \cite{huangshuo},  
 a finite difference scheme combined with Monte Carlo estimation has been
 introduced for fully nonlinear PDEs with gradients of order up to $2$
 using integration by parts.
 
 \medskip

 On the other hand,
 numerical solutions of semilinear PDEs have also
 been obtained by the multilevel Picard method
 \cite{hutzenthaler-mlp0}, 
 \cite{hutzenthaler-mlp2},
 \cite{hutzenthaler-mlp3},
 \cite{hutzenthaler-mlp1},
 with numerical experiments provided in \cite{hutzenthaler-mlp4}, 
 see also \cite{neufeld-wu-mlp} for a treatment of nonlocal PDEs.
 However, this approach is currently restricted to first order gradient
 nonlinearities. 
   
 \medskip

 Extending those techniques to nonlinearities in higher order derivatives
 involves several technical difficulties.
 In the case of branching diffusion approaches, 
 this involves a lack of integrability of the Malliavin-type weights used in 
 repeated integration by parts argument, 
 see page~199 of \cite{labordere}. This problem was also noted
 when dealing with pseudo-differential operators of the form
 $-\eta(\partial_x^2 )$ for the treatment of nonlocal PDEs \cite{penent}.
 
\medskip

In this paper, our method to deal with fully nonlinear PDEs of the form 
 \eqref{eq:1} relies on a marked branching process
 called a coding tree, represented by a random tree whose branches bear 
 operators, called \textit{codes}, instead of function values. 
 A general multiplicative functional 
 whose expected value provides a probabilistic representation
 of the PDE solution is then associated to the coding tree.
 To ensure the validity of the probabilistic representation we
 derive sufficient conditions on $f$ and $\phi$
 that ensure the finiteness of expected values. 

\medskip

Other probabilistic methods that can deal with
higher order derivatives usually involve pseudo-processes
created as limits of discrete random walks, see \cite{bonaccorsi}. 
However, the method developed in this paper is different,
as instead of creating a specific process whose generator
is behaving as a higher order derivative,
we use codes that carry information on the branches along
the tree. Once the tree leaves are reached, 
we make use of the code on the known terminal condition $\phi$
of the solution $u$. 

  \medskip 

  The idea of carrying information on nonlinearities along trees is not new, and
  has been developed in the case of ordinary differential equations (ODEs)
  in   \cite{butcher1963},
  see also, e.g.,  
  \cite{butcher2010}, 
  Chapters~4-6 of \cite{deuflhard}, 
  \cite{mclachlan}. 
  Butcher trees have found applications ranging from 
  geometric numerical integration to 
 stochastic differential equations, see for instance \cite{ehairer}
 and references therein, and 
 \cite{gubinelli}, \cite{bruned}, \cite{fossy}, for the use of
 decorated trees for stochastic partial differential equations, 
 and for their connections with the Butcher-Connes-Kreimer Hopf algebra
 \cite{connes}.

 \medskip
 
 In the approach of \cite{butcher1963}, the general idea is to write a Taylor expansion for the solution of a
 differential equation, and to represent every term using a specific tree structure.
 In this case, numerical evaluation of the solution
 requires to truncate the series by selecting certain trees.
  On the other hand, the stochastic branching method does not rely on truncations
  and can be used to estimate an infinite series as an
  expected value over almost surely finite random trees.  
  This approach has been applied in \cite{penent4} to the numerical
  estimation of ODE solutions 
  by the Monte Carlo method without the use of diffusion processes. 
  On the other hand, PDEs can be treated 
  by this method by attaching a random Brownian evolution
  to each tree branch. 

  \medskip 

In this paper, we provide probabilistic representations for
the solutions of a class of fully nonlinear parabolic PDEs of the form
\eqref{eq:1} 
 with functional nonlinearity $f(z_0,\ldots , z_n)$ in the solution $u$
 and its derivatives $\partial^k_x u(t,x)$, $k=1,\ldots , n$. 
 In the sequel, we
   denote by 
   \begin{equation}
\nonumber 
     \varphi (t,x) := \frac{\re^{- x^2 / ( 2t) }}{\sqrt{2\pi t} },
     \qquad x\in \real, 
     \end{equation} 
   the standard Gaussian kernel with variance $t>0$.
    We denote by ${\cal C}^k(\real^m)$ the set of $k$-times
 differentiable functions with continuous derivatives of orders up to
 $k \in \{0,\ldots , \infty\}$ on $\real^m$, and 
 for any $h\in {\cal C}^k(\real^m)$ and $(\lambda_1,\ldots , \lambda_n) \in \{0,\ldots , k\}^m$
  we use the notation 
 $$
 \partial_{z_1}^{\lambda_1}
\cdots
\partial_{z_m}^{\lambda_m} 
h(z_1,\ldots, z_m) 
 := \frac{\partial^{\lambda_1} }{\partial z_1^{\lambda_1}} \cdots \frac{\partial^{\lambda_m} }{\partial z_m^{\lambda_m}} h(z_1,\ldots , z_m), \quad m \geq 1. 
$$
Similarly, we denote by ${\cal C}^{1,k}([0,T]\times \real)$
functions $u(t,x)$ which are differentiable in time $t\in [0,T]$
and $k$ times differentiable in $x\in \real$ with continuous
partial derivatives, $1\leq k \leq \infty$.

 \noindent
    {\bf Assumption (\hypertarget{BGJhyp}{$A$})}.  
{\em
  Assume that
  \begin{enumerate}[i)]
    \item $f\in {\cal C}^\infty ( \real^{n+1})$ and $\phi \in {\cal C}^\infty ( \real )$,
 \item the PDE \eqref{eq:1} admits a unique solution $u\in {\cal C}^{1,\infty}([0,T]\times \real )$, written in integral or Duhamel formulation as 
\begin{eqnarray} 
 \label{Eint3-0}
 u (t,x)
 & = & \int_{-\infty}^\infty \varphi (T-t,y-x)\phi(y)dy
    \\
    \nonumber 
     & & +
   \int_t^T \int_{-\infty}^\infty \varphi (s-t,y-x)
     f\left(u(s,y) , \partial_y u (s,y), \ldots ,
  \partial^n_y u (s,y) \right)
  dy ds, 
\end{eqnarray} 
\item
  $\phi^{(k)}(u)\in \cap_{p=1}^{n+1} L^p(\real , \varphi ( \eta , x ) dx )$, $k\geq 0$, and
    $\partial^{\lambda_0}_{z_0} \cdots \partial^{\lambda_n}_{z_n} f(u,\partial_x u,\ldots , \partial_x^nu)
  \in \cap_{p=1}^{n+1} L^p([0,T]\times \real, \varphi (\eta ,x) dx ds)$, 
    $(\lambda_0,\ldots , \lambda_n) \in \inte^{n+1}$, 
for all $\eta > 0$. 
\end{enumerate} 
}
\noindent 
 We refer to e.g. Theorem~1.1 in \cite{krylov1983}
 for sufficient conditions for existence and uniqueness of
 smooth solutions to such fully nonlinear PDEs in the second order case.

 \medskip

   Starting from \eqref{Eint3-0}, 
   we will construct a random coding tree $\mathcal{T}_{t,x,c}$ 
   rooted at $(t,x)$ which is a random branching process
  driven by a standard Brownian motion $(W_t)_{t\in \real_+}$, 
 with branches bearing operators called codes and indexed by a set $\eufrak{C}$,
   such that the first branch of this tree bears the code $c = {\rm Id}$. 
   
   \medskip

   Next, we will construct a universal multiplicative functional $\mathcal{H}_\phi$ 
   of $\mathcal{T}_{t,x,c}$, such that the expectation 
   $u_c(t,x):=\E [\mathcal{H}_\phi(\mathcal{T}_{t,x,c}) ]$ solves the
  system of equations  
\begin{equation} 
\label{s1-1-0}
\left\{
\begin{array}{l}
  \displaystyle \partial_t u_c(t,x) + \frac{1}{2}\partial_x^2 u_c(t,x)
+
\sum_{Z \in \mathcal{M}(c)}
\prod_{z \in Z} u_z (t,x) = 0, \quad c\in \eufrak{C},  
\medskip
\\
 u_c (T,x) = c(u)(T,x), \qquad (t,x)\in [0,T] \times \real, 
\end{array}
\right.
\end{equation} 
 where $\mathcal{M}$ is a mapping called the \textit{mechanism}, 
 which sends any code $c\in \eufrak{C}$
 to a family of code tuples $Z \in {\cal M}(c)$ 
 which are associated to the new branches created in
 the random coding tree $\mathcal{T}_{t,x,c}$. 

 \medskip

 In Theorem~\ref{t1}, supposing in addition to Assumption~(\hyperlink{BGJhyp}{$A$}) 
 that the solution of the system \eqref{s1-1-0} is unique, 
 and given $T>0$ such that the 
 functional $\mathcal{H}(\mathcal{T}_{t,x,c})$ is integrable for all
 $(t,x) \in [0,T] \times \real$, 
   we derive a probabilistic representation of
   the form
\begin{equation}
\nonumber 
 u (t,x) :=  \E \big[ \mathcal{H}_\phi (\mathcal{T}_{t,x,{\rm Id}})\big],
 \quad (t,x) \in [0,T]\times \real,
\end{equation} 
for the solution $u(t,x)$ of 
\eqref{eq:1}. 
Sufficient conditions on $f$, $\phi$ for the boundedness of the functional $\mathcal{H}(\mathcal{T}_{t,x,c} )$ are derived Proposition~\ref{t3.3} under additional conditions on the probability density function $\rho$ of interbranching times in the random tree $\mathcal{T}_{t,x,c}$,
 over a sufficiently small time interval $[0,T]$. 
 
 \medskip

 In Section~\ref{s5} we present a Monte Carlo implementation of our algorithm
 for the numerical solutions of fully nonlinear PDEs on a sufficiently
 small time interval. 
 Numerical applications are presented 
 to semilinear, quasilinear and fully nonlinear PDEs.
 This includes in particular functional nonlinearities 
 which are not covered by standard branching methods
 that are designed for polynomial nonlinearities. 
 We also deal with examples involving higher order derivatives that 
 may not be treated by Malliavin-type integration by parts arguments 
 due to integrability issues, see page~199 of \cite{labordere},
 and are also not covered by 
 multilevel Picard methods, see e.g.
 \cite{hutzenthaler-mlp4},
 or BSDE methods, see e.g. 
 \cite{han2018solving},
 which are limited to first and second order gradients,
 respectively.

\medskip

Although our results are only valid in small time,
the numerical experiments performed in Section~\ref{s5}
for the Allen-Cahn equation \eqref{gl0}-\eqref{gl} 
and for the HJB equation \eqref{hjbgl0}, 
see     Tables~\ref{fjkldsf1}, \ref{t3}, \ref{t3-0}, \ref{fjl33}
    and
    Figures~\ref{ac} and \ref{f0-mlp}-$b)$, 
    show that the performance of our coding tree method compares
    favorably to those of
    the BSDE, branching diffusion, and MLP methods.
    In addition, some of our fully nonlinear
    examples, see Examples~3-a) and 3-b), are currently out of reach by other methods. 

\medskip

 This paper is organized as follows.
 Sections~\ref{s2} and \ref{s2.2} present the constructions of codes,
 mechanisms, and random coding trees. 
 In Section~\ref{s3} we state our main result
 Theorem~\ref{t1}
 which gives the probabilistic representation of the solution and
 its partial derivatives
 and give a sufficient condition that 
 ensures the integrability needed for the
 probabilistic representation of Theorem~\ref{t1} to hold.
 In Section~\ref{s5},
 we present numerical simulations that illustrate the method on specific examples.

\medskip

 The appendix contains a Mathematica implementation of the algorithm of Theorem~\ref{t1}
 in dimension one. 
 The Python codes designed for other numerical experiments are available
 at \url{https://github.com/nguwijy/coding_trees}. 
 
\subsubsection*{Preliminaries}
 For simplicity of exposition, Sections~\ref{s2}-\ref{s3}
 are presented in the one-dimensional
 case of PDEs of a single space variable $x\in \real$,
 while the codes used in Section~\ref{s5}
 are implemented in the $d$-dimensional setting.
 In the sequel we will use the following version of the multivariate
 Fa\`a di Bruno formula, which follows from Theorem~2.1 in \cite{constantine}.
\begin{prop}
  Let $n\geq 0$ and $k \geq 1$.
  Given $g\in {\cal C}^k(\real^{n+1})$ function of
  $(z_0,\ldots, z_n)$ 
  and $v\in {\cal C}^n([0,T]\times \real )$ function of
  $(t,x)$, we have 
\begin{align} 
  \nonumber
 &    \partial_x^k \big( g\big(v(t,x),\dots,\partial_x^n v(t,x)\big) \big) 
   \\
  \label{fdb}
  & = 
  k!
  \sum_{1 \leq \lambda_0 + \cdots + \lambda_n \leq k
    \atop 1 \leq s \leq k
  } 
\partial_{z_0}^{\lambda_0}
\cdots
\partial_{z_n}^{\lambda_n} 
g\big(v(t,x),\dots,\partial_x^n v(t,x)\big) 
\sum_{ 1 \leq | {\bf k}_1|,\ldots , |{\bf k}_s|, \ 1 \leq l_1 < \cdots < l_s
   \atop
          {
            k_1^i+\cdots + k_s^i = \lambda_i, \ 0 \leq i \leq n
            \atop
            |{\bf k}_1|l_1+\cdots + |{\bf k}_s|l_s = k
      }
  }
  \prod_{1 \leq j \leq s \atop 0 \leq q \leq n} 
  \frac{
   \big( \partial_x^{q+l_j} v(t,x) \big)^{k_j^q} 
    }{k_j^q!(l_j!)^{k_j^q}}
  , \qquad 
\end{align} 
 with 
${\bf k}_j := (k_j^0,\ldots , k_j^n)$, 
 $|{\bf k}_j|:= k_j^0+\cdots + k_j^n$
 and
 ${\bf k}_j!:= k_j^0!\cdots k_j^n!$,
 $j=1,\ldots , k$.
\end{prop} 
We will also need the Duhamel formula, which shows that 
the solution $v(t,x)$ of an equation of the form 
\begin{equation}
\nonumber 
\begin{cases}
  \displaystyle
  \partial_t v (t,x)
 + \frac{1}{2} \partial_x^2 v (t,x)
  + 
  g(t,x) 
  = 0,
  \medskip
  \\
v(T,x) = \phi(x), \quad (t,x) \in [0,T]\times \real, 
\end{cases}
\end{equation}
 can be represented in integral form as
\begin{equation}
  \label{duhamel} 
 v(t,x) = 
  \int_{-\infty}^\infty \varphi (T-t,y-x)
    \phi (y ) dy + \int_t^T \int_{-\infty}^\infty \varphi (s-t,y-x)
 g(s,y )dy ds,
\end{equation}
 $(t,x) \in [0,T]\times \real$.
\section{Codes and mechanism}
\label{s2}
 In this section we start by constructing the set of codes $\eufrak{C}$ based on the
 Duhamel formula \eqref{duhamel}, 
 and by iterations of the Duhamel formula we deduce the mechanism $\mathcal{M}$,
 which will results in the construction of a random coding tree $\mathcal{T}_{t,x,c}$
 on which every particle evolves according
 to a Brownian motion with generator is $\partial_x^2 /2$.

 \medskip

 In order to derive a probabilistic representation for the solution of
 \eqref{Eint3-0}, we will derive an integral formulation for
    $f\big(u(t,x),\partial_x u(t,x),\dots, \partial_x^n u(t,x)\big)$
    and iterate this process.
    In the sequel, given $h\in {\cal C}^\infty (\real^{n+1} )$
     we let $h^*$ denote the mapping 
 \begin{align}
\nonumber 
h^* : {\cal C}^{0,\infty} ([0,T]\times \real ) 
& \longrightarrow {\cal C}^{0,\infty} ([0,T]\times \real ) 
    \\
  \label{fdhskjfd}
  \psi = \{ (t,x) \mapsto \psi (t,x) \} & \longmapsto h^*(\psi ):= \big\{ (t,x) \mapsto h\big(\psi (t,x),\partial_x\psi (t,x), \ldots , \partial_x^n\psi (t,x)\big) \big\}, 
\end{align} 
 where $\real^{[0,T]\times \real}$ represents the set of functions from $[0,T]\times \real$ to $\real$,
 and for $k\geq 1$ we identify $\partial_x^k$ to the operator defined as 
\begin{align}
\nonumber 
    \partial_x^k : {\cal C}^{0,k}([0,T]\times \real ) & \longrightarrow {\cal C}^{0,0} ([0,T]\times \real )
    \\
\nonumber 
  \psi = \{ (t,x) \mapsto \psi (t,x) \} & \longmapsto \partial_x^k (\psi ):=\left\{
    (t,x) \mapsto
    \partial^k_x \psi (t,x) \right\}.  
\end{align}
Letting $v(t,x) := g\big(u(t,x),\partial_x u(t,x),\dots, \partial_x^n u(t,x)\big)$ where
$g\in {\cal C}^\infty (\real^{n+1})$, 
    by the Fa\`a di Bruno formula \eqref{fdb} 
    we have 
\begin{align} 
  \nonumber
  & \partial_tv(t,x)
+ \frac{1}{2}\partial_x^2 v(t,x)
= \sum_{k=0}^n \partial_{z_k} g\big(u(t,x),\dots,\partial_x^n u(t,x)\big)
\partial_x^k \left( \partial_t u(t,x) + \frac{1}{2}\partial_x^2 u(t,x) \right)
\\ 
\nonumber
   &  \quad
+ \frac{1}{2}\sum_{k=0}^n \sum_{l=0}^n \partial_{z_l} \partial_{z_k} g\big(u(t,x),\dots,\partial_x^n u(t,x)\big) \partial_x^{k+1} u(t,x) \partial_x^{l+1} u(t,x)
\\
\nonumber
&=  - \sum_{k=0}^n \partial_{z_k} g\big(u(t,x),\dots,\partial_x^n u(t,x)\big) \partial_x^k f\big(u(t,x),\dots,\partial_x^n u(t,x)\big)
\\
\nonumber
  &   \quad  +  \frac{1}{2}\sum_{k=0}^n \sum_{l=0}^n \partial_{z_l} \partial_{z_k} g\big(u(t,x),\dots,\partial_x^n u(t,x)\big) \partial_x^{k+1} u(t,x) \partial_x^{l+1} u(t,x)
\\
\nonumber
    &= - \partial_{z_0} g\big(u(t,x),\dots,\partial_x^n u(t,x)\big)
  f\big(u(t,x),\dots,\partial_x^n u(t,x)\big)
  \\
\nonumber
    & \quad - \sum_{k=1}^n k! \partial_{z_k} g\big(u(t,x),\dots,\partial_x^n u(t,x)\big)
\hskip-0.2cm
\sum_{1\leq |\lambda| \leq k
    \atop 1 \leq s \leq k
  }
    \big( \partial_{z_0}^{\lambda_0} \cdots \partial_{z_n}^{\lambda_n} f\big)^*
  \hskip-1cm
  \sum_{ 1 \leq | {\bf k}_1|,\ldots , |{\bf k}_s| , \ 
     1 \leq l_1 < \cdots < l_s
      \atop
          {
            k_1^i+\cdots + k_s^i = \lambda_i, \ 0 \leq i \leq n
            \atop
            |{\bf k}_1|l_1+\cdots + |{\bf k}_s|l_s = k
            }
    }
   \prod_{1 \leq j \leq s\atop
    0 \leq q \leq n}
  \frac{
  \big( \partial_x^{q+l_j} u(t,x)\big)^{k_j^q} }{k_j^q!(l_j!)^{k_j^q}}
  \\
\label{fdjskf1} 
    & \quad 
   +  \frac{1}{2}\sum_{j=0}^n \sum_{l=0}^n \partial_{z_l} \partial_{z_j} g\big(u(t,x),\dots,\partial_x^n u(t,x)\big) \partial_x^{j+1} u(t,x) \partial_x^{l+1} u(t,x), 
\end{align} 
where $|\lambda|:=\lambda_0+\cdots + \lambda_n$,  $\lambda \in \inte^{n+1}$. 
The PDE \eqref{fdjskf1} 
 can be rewritten in integral form by the Duhamel formula \eqref{duhamel} as 
\begin{align} 
\label{kdsf}
& 
     g^*(u)(t,x) = g\big(u(t,x),\partial_x u(t,x),\dots, \partial_x^n u(t,x)\big)
      \\
      \nonumber
      & 
  = 
  \int_{-\infty}^\infty \varphi (T-t,y-x)
    g\big(\phi (y),\partial_y \phi(y),\dots, \partial_y^n \phi(y)\big)
  dy
    \\
  \nonumber
   & \quad 
 + \int_t^T \int_{-\infty}^\infty \varphi (s-t,y-x)
 \Bigg(
 \partial_{z_0} g\big(u(s,y),\dots,\partial_y^n u(s,y)\big)
  f\big(u(s,y),\dots,\partial_y^n u(s,y)\big)
  \\
  \nonumber
  & \quad \quad
  \hskip-0.1cm
  + \sum_{k=1}^n k! \partial_{z_k} g\big(u(s,y),\dots,\partial_y^n u(s,y)\big)
  \hskip-0.1cm
  \sum_{1\leq |\lambda| \leq k
    \atop 1 \leq s \leq k
  }
  \hskip-0.1cm
  \big( \partial_{z_0}^{\lambda_0} \cdots \partial_{z_n}^{\lambda_n} f\big)^*
  \hskip-0.7cm
  \sum_{ 1 \leq | {\bf k}_1|,\ldots , |{\bf k}_s|, \ 
     1 \leq l_1 < \cdots < l_s
      \atop
          {
            k_1^i+\cdots + k_s^i = \lambda_i, \ 0 \leq i \leq n
            \atop
            |{\bf k}_1|l_1+\cdots + |{\bf k}_s|l_s = k
            }
      } 
\hskip-0.1cm
     \prod_{1 \leq j \leq s
    \atop
  0 \leq q \leq n}
\hskip-0.1cm
  \frac{ 
  \big( \partial_y^{q+l_j} u(s,y) \big)^{k_j^q} 
}{k_j^q!(l_j!)^{k_j^q}}
  \\
  \nonumber
  & \quad  \quad \quad  \left.
   - \frac{1}{2}\sum_{j=0}^n \sum_{l=0}^n \partial_{z_l} \partial_{z_j} g\big(u(s,y),\dots,\partial_y^n u(s,y)\big) \partial_y^{j+1} u(s,y) \partial_y^{l+1} u(s,y)\right) dy ds.
\end{align} 
 In order to formalize and extend the above iteration 
 we introduce the following definition, which relies on
 \eqref{fdhskjfd}. 
\begin{definition}
 We let $\eufrak{C}$ denote the set of 
 operators from ${\cal C}^{0,\infty} ([0,T]\times \real )$ to
 ${\cal C}^{0,\infty} ([0,T]\times \real )$, 
 called \textit{codes}, and defined as  
$$
 \eufrak{C}:= \left\{
         {\rm Id}, \ 
         \
         \big( a \partial_{z_0}^{\lambda_0} \cdots \partial_{z_n}^{\lambda_n} f\big)^*, 
                 \ \partial_x^k \ : \ \lambda 
                 \in \inte^{n+1} ,
                 \ a \in \real\setminus \{0\},
                                  \ k \geq 1 
                          \right\}, 
$$
         where
                  ${\rm Id}$
         denotes the identity on ${\cal C}^{0,\infty} ([0,T]\times \real )$.
\end{definition}
The role of the parameter $a\in \real\setminus \{0\}$ appearing in the
  definition of $\eufrak{C}$ is to account for possible real
  coefficients appearing in front of partial derivatives in
  the mapping $\mathcal{M}$, called the \textit{mechanism}, 
  defined on $\eufrak{C}$ according to \eqref{fdb} and \eqref{kdsf},
  by matching a \textit{code} $c\in \eufrak{C}$
 to a set $\mathcal{M}(c)$ of code tuples.
\begin{definition} 
 The mechanism $\mathcal{M}$ is defined on $\eufrak{C}$ by 
 letting $\mathcal{M} ( {\rm Id} ) := \{ f^* \}$, and 
\begin{align} 
\label{fjkdsl}
 \mathcal{M} ( g^* ) & :=  \big\{ 
   (f^*,(\partial_{z_0} g)^*)
 \big\}
  \\
  \nonumber 
&
  \qquad
 \bigcup
 \bigcup_{
  1 \leq \lambda_0 + \cdots + \lambda_n \leq k 
    \atop
 { 1 \leq | {\bf k}_1|,\ldots , |{\bf k}_s| , 
 1 \leq l_1 < \cdots < l_s
      \atop
          {
            k_1^i+\cdots + k_s^i = \lambda_i, \ 0 \leq i \leq n
            \atop
            |{\bf k}_1|l_1+\cdots + |{\bf k}_s|l_s = k, 
            \ 1 \leq s \leq k \leq n }
                }
}
   \left\{            \left(
       (\partial_{z_k} g)^*,
  k!
   \big( \partial_{z_0}^{\lambda_0} \cdots \partial_{z_n}^{\lambda_n} f\big)^*
           , 
   \left(
   \frac{\partial_x^{q+l_j}}{
   k_j^q!(l_j!)^{k_j^q}
   }
   \right)_{r=1,\ldots , k_j^q
   \atop
   {j=1,\ldots , s \atop 
     q= 0, \ldots n }
 } 
  \right)
  \right\}
  \\
  \nonumber 
&
  \qquad
   \bigcup
   \bigcup_{j,l=0 ,\ldots , n}
  \left\{ 
  \left(
  - \frac{1}{2}(\partial_{z_l} \partial_{z_j} g)^* , \partial_x^{j+1} ,  \partial_x^{l+1} 
  \right)  \right\}
  , \qquad g^*\in \eufrak{C},
\end{align} 
 and
\begin{eqnarray} 
  \nonumber
  \mathcal{M} \big( \partial_x^k \big) :=  
  \bigcup_{
  1 \leq \lambda_0 + \cdots + \lambda_n \leq k 
  , \ s=1,\ldots , k
  \atop
{
  1 \leq | {\bf k}_1|,\ldots , |{\bf k}_s| , \ 
 1 \leq l_1 < \cdots < l_s 
      \atop
          {
            k_1^i+\cdots + k_s^i = \lambda_i, \ 0 \leq i \leq n
            \atop
            |{\bf k}_1|l_1+\cdots + |{\bf k}_s|l_s = k 
            }
      }
}
  \left\{
\left(
    k!
   \big( \partial_{z_0}^{\lambda_0} \cdots \partial_{z_n}^{\lambda_n} f\big)^*
           , 
   \left(
   \frac{\partial_x^{q+l_j}}{
   k_j^q!(l_j!)^{k_j^q}
   }
   \right)_{r=1,\ldots , k_j^q
   \atop
  {j=1,\ldots , s \atop 
  q= 0, \ldots n }
   }
   \right) 
      \right\},
   \quad k \geq 1. 
   \\
   \label{fjkdsl-8}
   \end{eqnarray} 
\end{definition} 
\subsubsection*{Example - semilinear PDEs} 
 As an example, let $n=0$ and consider a semilinear PDE of the form 
\begin{equation}
   \label{eq:112}
\begin{cases}
  \displaystyle
  \partial_t u(t,x) + \frac{1}{2}\partial_x^2 u(t,x) + f(u(t,x)) = 0
  \medskip
  \\
u(T,x) = \phi(x), \qquad (t,x) \in [0,T] \times \real. 
\end{cases}
\end{equation} 
 Letting $v(t,x) := g(u(t,x))$, Equation~\eqref{fdjskf1} reads 
\begin{equation*}
\begin{split}
  \partial_t v (t,x)+ \frac{1}{2}\partial_x^2 v (t,x)
  &=
  \left( \partial_t u (t,x)+ \frac{1}{2} \partial_x^2 u(t,x) \right)
  g'(u(t,x)) 
  +  \frac{1}{2} (\partial_x u(t,x))^2 g''(u(t,x))
  \\
&= -f(u(t,x)) g' (u(t,x)) + \frac{1}{2}(\partial_x u(t,x))^2 g''(u(t,x)). 
\end{split}
\end{equation*}
 Therefore, by Duhamel's formula \eqref{duhamel},
 $v(t,x)$ satisfies  the integral equation 
 \begin{align*}
   v(t,x) & = \int_{-\infty}^\infty \varphi (T-t,y-x) g( \phi(y) ) dy
   \\
   & \quad 
   + \int_t^T \int_{-\infty}^\infty \varphi (s-t,y-x)
 \left( f(u(s,y)) g' (u(s,y)) - \frac{1}{2}(\partial_y u(s,y))^2 g''(u(s,y)) \right) dy ds 
\end{align*} 
 as in \eqref{kdsf},
 and the set of $\eufrak{C}$ codes is given by 
$$
 \eufrak{C} := \big\{
           {\rm Id}, \ \partial_x, \ a f^{(k)}, \ a \in \real\setminus \{0\}, \ k \in \mathbb{N} \big\}. 
$$
                 In this case, the mechanism $\mathcal{M}$ is given by 
                 \begin{equation}
                   \label{jfklds} 
 {\cal M}({\rm Id}) := \{f^* \},
 \ \ 
     {\cal M}(g^*) := \left\{\big(f^*,( g' )^*\big)
     ;
      \left(\partial_x,\partial_x, -\frac{1}{2}( g'' )^* \right)\right\}, 
 \ \ 
 {\cal M}(  \partial_x ) := \big\{ \big( (f' )^* ,\partial_x \big) \big\},  
\end{equation} 
                 for $g\in {\cal C}^\infty ( \real^{n+1})$
                  of the form $g = a f^{(k)}$, $a \in \real\setminus \{0\}$, $k \geq 0$.
 In this case, every code tuple in $\mathcal{M}(c)$ has at most $3$ elements for any $c\in \eufrak{C}$, and the time complexity of the algorithm can be estimated from the mean depth of the random tree $\mathcal{T}_{0,x,c}$, which grows exponentially as a function of $T>0$, see e.g. \S~4 of \cite{penent4}. 
\subsubsection*{Example - first order gradient nonlinearity} 
 As a second example, let $n=1$ and consider the nonlinear PDE 
\begin{equation}
\nonumber 
\begin{cases}
  \displaystyle \partial_t u(t,x) + \frac{1}{2}\partial_x^2 u(t,x) + f(u(t,x),\partial_x u(t,x)) = 0
  \medskip
  \\
u(T,x) = \phi(x). 
\end{cases}
\end{equation} 
 Letting $v(t,x) := g(u(t,x),\partial_x u(t,x))$, Equation~\eqref{fdjskf1} reads 
\begin{eqnarray*}
\lefteqn{ 
   \partial_t v(t,x)
  +\frac{1}{2} \partial_x^2 v(t,x)
  = \partial_{z_0} g(u(t,x),\partial_x u(t,x)) 
  \left( \partial_t u(t,x) + \frac{1}{2}\partial_x^2 u(t,x) \right) 
}
\\
 & & + \partial_{z_1} g(u(t,x),\partial_x u(t,x))  \left(
\partial_x \partial_t u(t,x) + \frac{1}{2} \partial_x^3 u(t,x) \right)
\\
& & + \frac{1}{2}\partial_{z_0}^2 g(u(t,x),\partial_x u(t,x)) (\partial_x u(t,x))^2+ \frac{1}{2} \partial_{z_1}^2 g(u(t,x),\partial_x u(t,x)) (\partial_x^2 u(t,x))^2 \\
& & + \partial_{z_0}\partial_{z_1} f(u(t,x),\partial_x u(t,x)) \partial_x u(t,x) \partial_x^2 u(t,x)\\
 &=& - \partial_{z_0} g(u(t,x),\partial_x u(t,x)) f(u(t,x),\partial_x u(t,x))\\
 && -  \partial_{z_1} g(u(t,x),\partial_x u(t,x)) \big( \partial_{z_0} f(u(t,x),\partial_x u(t,x)) \partial_x u(t,x) + \partial_{z_1} f(u(t,x),\partial_x u(t,x)) \partial_x^2 u(t,x) \big)
\\
&  &+ \frac{1}{2} \partial_{z_0}^2 g(u(t,x),\partial_x u(t,x))
(\partial_x u(t,x))^2+ \frac{1}{2}\partial_{z_1}^2 g(u(t,x),\partial_x u(t,x)) (\partial_x^2 u(t,x))^2 \\
& & + \partial_{z_0}\partial_{z_1} g(u(t,x),\partial_x u(t,x)) \partial_x u(t,x) \partial_x^2 u(t,x),
\end{eqnarray*} 
 and by Duhamel's formula \eqref{duhamel},
 $v(t,x)$ satisfies  the integral equation 
 \begin{align}
\nonumber
   v(t,x) & = \int_{-\infty}^\infty \varphi (T-t,y-x)
   g(u(t,y),\partial_y u(t,y)) dy
   \\
   \nonumber
   & 
      + \int_t^T \int_{-\infty}^\infty \varphi (s-t,y-x)
      \Big(
      \partial_{z_0} g(u(t,y),\partial_y u(t,y)) f(u(t,y),\partial_y u(t,y))
      \\
\nonumber
    & + \partial_{z_1} g(u(t,y),\partial_y u(t,y)) \big( \partial_{z_0} f(u(t,y),\partial_y u(t,y)) \partial_y u(t,y) + \partial_{z_1} f(u(t,y),\partial_y u(t,y)) \partial^2_y u(t,y) \big)
\\
\nonumber
   &   - \frac{1}{2} \partial_{z_0}^2 g(u(t,y),\partial_y u(t,y))
(\partial_y u(t,y))^2
 - \frac{1}{2}\partial_{z_1}^2 g(u(t,y),\partial_y u(t,y)) \big(\partial^2_y u(t,y) \big)^2 \\
   \label{fjdksl342} 
   &  - \partial_{z_0}\partial_{z_1} g(u(t,y),\partial_y u(t,y)) \partial_y u(t,y) \partial^2_y u(t,y)
\Big) dy ds 
\end{align} 
as in \eqref{kdsf}.
 In this case, the set of codes is given by 
$$
\eufrak{C} := \big\{ {\rm Id}, \ \partial_x^k , \ ( a \partial_{z_0}^l \partial_{z_1}^m f)^* \ : \ a \in \real\setminus \{0\}, \ k \geq 1, \ l, m \geq 0 \big\}, 
$$
  and 
  by \eqref{fjdksl342}, the mechanism $\mathcal{M}$ satisfies 
  ${\cal M}({\rm Id}) := \{f^*\}$, 
  \begin{align*} 
    & {\cal M}(g^*) := \big\{ \big(f^*, (\partial_{z_0} g )^*\big)
    ;
    \big( (\partial_{z_1} g)^*, ( \partial_{z_0} f)^*, \partial_x \big)
    ;
    \big( (\partial_{z_1} g)^*, (\partial_{z_1} f)^*, \partial_x^2 \big)
    ; 
 \\
  &
   \ 
   \left. 
  \left(-\frac{1}{2} ( \partial_{z_0}^2 g )^*,\partial_x,\partial_x\right) ;
  \left(-\frac{1}{2} ( \partial_{z_1}^2 g)^*,\partial_x,\partial_x\right) ; 
  \left(-\frac{1}{2} (\partial_{z_0}\partial_{z_1}g)^*,\partial_x,\partial_x^2 \right)
  ;
  \left(-\frac{1}{2} (\partial_{z_0}\partial_{z_1}g)^*,\partial_x,\partial_x^2 \right)
\right\} 
\end{align*} 
for $g\in {\cal C}^\infty ( \real^{n+1})$ of the form $g = a \partial_{z_0}^k \partial_{z_1}^l f$, $a \in \real\setminus \{0\}$, $k,l \geq 0$, and 
$$
 {\cal M}(\partial_x ) := \left\{ \big( ( \partial_{z_0} f )^*, \partial_x \big)
  ;
 \big( (\partial_{z_1} f )^*, \partial_x^2 \big)  \right\}.
 $$ 
 This makes it possible to find the image of $\partial_x^k$ as well,
  for example we have 
\begin{align*} 
  {\cal M}( \partial_x^2 ) & =
  \big\{ \big( (\partial_{z_0} f )^*, \partial_x^2 \big);
  \big( (\partial_{z_0}^2 f )^*, \partial_x, \partial_x \big);
  \big( (\partial_{z_0}\partial_{z_1} f)^*, \partial_x^2, \partial_x \big);
  \big( (\partial_{z_1} f)^*, \partial_x^3 \big);
  \big( ( \partial_{z_0}\partial_{z_1} f)^*,\partial_x,\partial_x^2 \big);
  \\
  &
  \qquad \big( ( \partial_{z_1}^2 f)^*,\partial_x^3 ,\partial_x^2 \big)  \big\}. 
  \end{align*} 
 We close this section with the following key lemma which
 shows that $c(u)$ satisfies a system
 of equations indexed by $c\in \eufrak{C}$,
 and present its application to a semilinear example.
 \begin{lemma}
  \label{l1}
  For any code $c\in \eufrak{C}$ we have
\begin{equation} 
 \label{s1}
 c(u)(t,x) = \int_{-\infty}^\infty \varphi (T-t,y-x) c ( u)(T,y) dy 
+ 
\sum_{Z \in \mathcal{M}(c)}
\int_t^T \int_{-\infty}^\infty \varphi (s-t,y-x)
\prod_{z \in Z}  z(u)(s,y) dy ds, 
\end{equation} 
$(t,x)\in [0,T]\times \real$.
 \end{lemma} 
 \begin{Proof}
  When $c={\rm Id}$ we have
\begin{align*} 
 & 
    {\rm Id}(u)(t,x)  =  u(t,x)
  \\
   & = 
  \int_{-\infty}^\infty \varphi (T-t,y-x) \phi (y) dy 
+ 
\int_t^T \int_{-\infty}^\infty \varphi (s-t,y-x)
  f\left(u(s,y) , \partial_y u (s,y), \ldots ,
  \partial^n_y u (s,y) \right)
 dy ds
  \\
   & = 
  \int_{-\infty}^\infty \varphi (T-t,y-x) u(T,y) dy 
+ 
\int_t^T \int_{-\infty}^\infty \varphi (s-t,y-x)
  f^* (u) (s,y) dy ds
  \\
   & = 
  \int_{-\infty}^\infty \varphi (T-t,y-x) {\rm Id}(u)(T,y) dy 
+ 
\int_t^T \int_{-\infty}^\infty \varphi (s-t,y-x)
  f^* (u) (s,y) dy ds
,
\end{align*}
 hence \eqref{s1} holds since $\mathcal{M} ( {\rm Id} ) = \{ f^* \}$. 
 When $c\not= {\rm Id}$ in written as $c = g^*\in \eufrak{C}$, the equation 
 \eqref{kdsf}
 satisfied by $g^*(y)(t)$ reads 
 \begin{eqnarray}
\nonumber 
  \lefteqn{
  g^*(u)(t,x)
  = 
  \int_{-\infty}^\infty \varphi (T-t,y-x)
  g^*(u)(T,y)
  dy
  }
  \\
  \nonumber
  & &
  + \int_t^T \int_{-\infty}^\infty \varphi (s-t,y-x)
 \Bigg(
 ( \partial_{z_0} g )^* (u) (s,y) f^* (u)(s,y)
  \\
  \nonumber
  & & \quad + \sum_{k=1}^n k! (\partial_{z_k} g)^*(u)(s,y)
  \sum_{1\leq |\lambda| \leq k
    \atop 1 \leq s \leq k
  }
    \big( \partial_{z_0}^{\lambda_0} \cdots \partial_{z_n}^{\lambda_n} f\big)^* (s,y)
  \sum_{ 1 \leq | {\bf k}_1|,\ldots , |{\bf k}_s|, \ 
     1 \leq l_1 < \cdots < l_s
     \atop
          {
            k_1^i+\cdots + k_s^i = \lambda_i, \ 0 \leq i \leq n
            \atop
            |{\bf k}_1|l_1+\cdots + |{\bf k}_s|l_s = k
            }
    }
  \prod_{1 \leq j \leq s
    \atop
 0 \leq q \leq n
  } 
  \frac{
   \big( \partial_x^{q+l_j} (u) (s,y) \big)^{k_j^q} 
  }{k_j^q!l_j!^{k_j^q}}
  \\
  \nonumber
  & & \left.
   \quad - \frac{1}{2}\sum_{j=0}^n \sum_{l=0}^n ( \partial_{z_l} \partial_{z_j} g)^*(u)(s,y) \partial_x^{j+1} (u) (s,y) \partial_x^{l+1} (u) (s,y)\right) dy ds. 
\end{eqnarray}
 Finally, when $c=\partial_x^k$, $k\geq 1$, the Fa\`a di Bruno formula shows that 
\begin{align}
\nonumber 
    & 
  \partial_x^k (u)(t,x) =
  \int_{-\infty}^\infty \varphi (T-t,y-x)
  \partial_x^k (u)(T,y)
  dy
    \\
  \nonumber
 & 
 + k! \int_t^T \int_{-\infty}^\infty \varphi (s-t,y-x)
 \hskip-0.1cm
 \sum_{1\leq |\lambda| \leq k
   \atop 1 \leq s \leq k
 } 
 \hskip-0.1cm
  \big( \partial_{z_0}^{\lambda_0} \cdots \partial_{z_n}^{\lambda_n} f\big)^* (s,y) 
 \hskip-0.7cm
 \sum_{ 1 \leq | {\bf k}_1|,\ldots , |{\bf k}_s|,
    \  1 \leq l_1 < \cdots < l_s
      \atop
          {
            k_1^i+\cdots + k_s^i = \lambda_i, \ 0 \leq i \leq n
            \atop
            |{\bf k}_1|l_1+\cdots + |{\bf k}_s|l_s = k
            }
    }
  \prod_{1 \leq j \leq s
    \atop
    0 \leq q \leq n
  } 
  \hskip-0.2cm
  \frac{ \big( \partial_x^{q+l_j} (u) (s,y) \big)^{k_j^q}
  }{k_j^q!l_j!^{k_j^q}}
 dy ds, 
\end{align} 
and \eqref{s1} follows from the definition \eqref{fjkdsl} of $\mathcal{M}$.
The exchange between summation over $Z\in \mathcal{M}(c)$
and integrals is justified by
Assumption~(\hyperlink{BGJhyp}{$A$})-$(iii)$. 
\end{Proof}

\subsubsection*{Example - semilinear PDEs} 
 In the case of a semilinear PDE of the form \eqref{eq:112} with $n=0$, 
 by the Duhamel formulation for the PDE satisfied by the functions $u(t,x)$,
 $\partial_x u(t,x)$, $a f^{(k)} (u(t,x))$,
 the system of equations \eqref{s1} reads 
  $$
  \left\{
  \begin{array}{l}
    \displaystyle
    u(t,x) = \int_{-\infty}^\infty \varphi (T-t,y-x) \phi(y) dy + \int_t^T \int_{-\infty}^\infty \varphi (s-t,y-x) f(u(s,y)) dy ds
    \medskip
    \\
    \displaystyle
  a f^{(k)} ( u(t,x) )= \int_{-\infty}^\infty \varphi (T-t,y-x) a f^{(k)} ( \phi (y) ) dy
\medskip
    \\
 \quad 
    \displaystyle
 + \int_t^T \int_{-\infty}^\infty \varphi (s-t,y-x)
\left( af ( u(s,y)) f^{(k+1)}  ( u(s,y)) - \frac{a}{2}(\partial_y u (s,y))^2 f^{(k+2)}  ( u (s,y)) \right) dy ds  
\medskip
    \\
    \displaystyle
\partial_x u(t,x) =  \int_{-\infty}^\infty \varphi (T-t,y-x) \partial_x \phi(y) dy + \int_t^T \int_{-\infty}^\infty \varphi (s-t,y-x) f' (  u (s,y) ) \partial_y u (s,y) dy  ds, 
  \end{array}
  \right. 
  $$
  $a \in \real\setminus \{0\}$, $k \in \mathbb{N}$,
  where the last equation can also be obtained by applying Duhamel's formula
  \eqref{Eint3-0} to
  $$
  \partial_t \partial_x u(t,x) + \frac{1}{2}\partial_x^2 \partial_x u(t,x)
  + f' (u(t,x)) \partial_x u(t,x) = 0.
  $$
 
\section{Random coding trees}
\label{s2.2}
This section introduces the random coding trees used for the
probabilistic representation of PDE solutions.
Let $\rho: \mathbb{R}^+ \rightarrow (0,\infty)$ be a
probability density function on $\real_+$.
For each $c \in \eufrak{C}$ we let $I_c$ be a random variable taking values uniformly in $\mathcal{M}(c)$ and note $q_c(b) := \mathbb{P}(I_c=b) $ where $b\in \mathcal{M}(c)$. 
 In addition, we consider 
\begin{itemize}
\item an i.i.d. family $(\tau^{i,j})_{i,j\geq 1}$ of random variables
 with probability density function $\rho$ on $\real_+$,
\item for each $c \in \eufrak{C}$,
  an independent family $(I_c^{i,j})_{i,j\geq 1}$ of i.i.d. discrete
  random variables on the finite set ${\cal M}(c)$, with distribution 
  $$
  \mathbb{P}\big( I_c^{i,j}=b \big) = q_c(b) >0,
  \qquad b \in \mathcal{M}(c), 
  $$
\item
  an independent family $(W^{i,j})_{i,j\geq 1}$
  of Brownian motions. 
\end{itemize}
 In addition, the sequences $(\tau^{i,j})_{i,j\geq 1}$, $(I_c^{i,j})_{c\in \eufrak{C}, i,j\geq 1}$ and $(W^{i,j})_{i,j\geq 1}$ are assumed to be mutually independent.

\medskip

We consider a coding branching process starting
from a particle $x\in \real$ at time $t \in [0,T]$ with label $\bar{1}=(1)$,
which evolves according to
the process $X_{s,x}^{\bar{1}} = x + W_{s-t}^{1,1}$, $s \in [t,t+\tau^{1,1}]$
and bears a code $c \in \eufrak{C}$.
If $\tau^{1,1}<T-t$, the process branches at time $t+\tau^{1,1}$
into new independent copies of
$(W_t)_{t \in \real_+}$, each of them started at
$X_{t+\tau^{1,1}}$ at time $t+\tau^{1,1}$.
Based on the value of $\big|I_c^{1,1}\big| \in \mathbb{N}$,
a family of $\big|I_c^{1,1}\big|$ new branches
are created. If $I_c^{1,1} = ( c_1,\ldots ,c_l )$ the $i$-th new branch will bear the code $c_i$, $i=1,\ldots , l$.

\medskip
 
Every new particle then follows independently
another copy of the same branching process as the initial particle, and
every branch stops when it reaches the horizon time $T$.
Particles at generation $n\geq 1$ are assigned a label of the form
 $\bar{k} = (1,k_2,\ldots ,k_n) \in \mathbb{N}^n$,
and their parent is labeled $\bar{k}- := (1,k_2,\ldots ,k_{n-1})$.
The particle labeled $\bar{k}$ is born at time $T_{\bar{k}-}$
and its lifetime $\tau^{n,\pi_n(\bar{k})}$ is the element of index
$\pi_n(\bar{k})$ in the i.i.d. sequence
$(\tau^{n,j})_{j\geq 1}$,
defining an injection
$$\pi_n:\mathbb{N}^n \to \mathbb{N},
\qquad n\geq 1.
$$ 
The random evolution of particle $\bar{k}$
is given by
\begin{equation}
\nonumber 
X_{s,x}^{\bar{k}} := X^{\bar{k}-}_{T_{\bar{k}-},x}+W_{s-T_{\bar{k}-}}^{n,\pi_n(\bar{k})},
\qquad s\in [T_{\bar{k}-},T_{\bar{k}}],
\end{equation} 
  where $T_{\bar{k}} := \T_{\bar{k}-} + \tau^{n,\pi_n(\bar{k})}$.

   \medskip 

 If $T_{\bar{k}} := \T_{\bar{k}-} + \tau^{n,\pi_n(\bar{k})} < T$,
 we draw a sample
 $I_c^{n,\pi_n(\bar{k})} = (c_1,\ldots ,c_l)$ uniformly in $\mathcal{M}(c)$,
 and the particle $\bar{k}$ branches into
 $\big|I_c^{n,\pi_n(\bar{k})}\big|$ offsprings at generation $(n+1)$,
 which are labeled by $\bar{k}=(1,\ldots ,k_n,i)$, $i=1,\ldots ,\big|I_c^{n,\pi_n(\bar{k})}\big|$.
 The particle with label ending with an integer $i$ will carry the code $c_i$,
 and the code of particle $\bar{k}$ is denoted by
 $c_{\bar{k}} \in \eufrak{C}$. 
 The labels are only used to distinguish the particles in the
 branching process. 

 \medskip 

 The set of particles dying before time $T$ is denoted by $\mathcal{K}^{\circ}$,
whereas those dying after $T$ form a set denoted by $\mathcal{K}^{\partial}$.
\begin{definition}
 When started at time $t\in [0,T]$ from
 a position $x\in \real$ and a code 
 $c \in \eufrak{C}$ on its first branch,
 the above construction yields
 a marked branching process called a random coding tree, and 
 denoted by $\mathcal{T}_{t,x,c}$.
\end{definition}
We note that 
the random branching tree $\mathcal{T}_{t,x,c}$ is non-explosive in
finite time since the number of branching times is a.s. finite, as
the sequence $(\tau^{i,j})_{i,j\geq 1}$ is i.i.d.. 
 The random tree $\mathcal{T}_{t,x,{\rm Id}}$
 will be used for the stochastic representation of the solution
 $u(t,x)$ of the PDE \eqref{eq:1}, while the trees $\mathcal{T}_{t,x,c}$ 
 will be used for the stochastic representation of $c (u) (t,x)$. 
 The next table summarizes the notation introduced so far.

\medskip

\begin{table}[H] 
  \centering
\footnotesize 
  \begin{tabular}{||l | c||}
 \hline
 Object & Notation \\ [0.5ex]
 \hline\hline
 Initial time  & $t$ \\
 \hline
 Initial position & $x$ \\
 \hline
 Tree rooted at $(t,x)$ with initial code $c$ &  $\mathcal{T}_{t,x,c}$ \\  \hline
 Particle (or label) of generation $n\geq 1$ & $\bar{k}=(1,k_2,\ldots ,k_n)$\\
 \hline
 First branching time & $T_{\bar{1}}$\\
 \hline
 Lifespan of a particle & $\tau_{\bar{k}} = T_{\bar{k}} - T_{\bar{k}-}$ \\
 \hline
Birth time of a particle $\bar{k}$ & $T_{\bar{k}-}$ \\
\hline
Death time of a particle $\bar{k}$ & $T_{\bar{k}}$ \\
\hline
Position at birth & $X^{\bar{k}}_{T_{\bar{k}-},x}$\\
\hline
Position at death & $X^{\bar{k}}_{T_{\bar{k}},x}$ \\
\hline
Code of a particle $\widebar{k}$ & $c_{\bar{k}}$
\\ 
\hline
\end{tabular}
\caption{Summary of branching tree notation.}
\end{table} 
\subsubsection*{Example - semilinear PDEs} 
 In the case of a semilinear PDE of the form \eqref{eq:112} with $n=0$, 
 the distributions 
  $q_c$, $c\in \eufrak{C}$, on the mechanism \eqref{jfklds} is
  given by 
 $q_{\rm Id} (f^*)=1$,
 $q_{\partial_x} ( ( (f')^* ,\partial_x ) ) =1$, and  
$$
  \left\{
  \begin{array}{l}
    \displaystyle
    q_{a f^{(k)} } \big( \big(f^*,\big(a f^{(k+1)}  \big)^*\big)\big) =
  \P\big( I_{a f^{(k)}  }
  =
  \big(f^*,\big(a f^{(k+1)}  \big)^*\big)\big) = \frac{1}{2},
  \\
  \\
 \displaystyle
    q_{a f^{(k)} } \left( 
 \left(\partial_x,\partial_x, -\frac{1}{2}(a f^{(k+2)}  )^* \right)
  \right)
  =
 \P \left( 
 I_{a f^{(k)} }
 =
 \left(\partial_x,\partial_x, -\frac{1}{2}(a f^{(k+2)}  )^* \right)
  \right)
=\frac{1}{2}, \quad k \geq 0.
  \end{array}
  \right.
$$ 

  \noindent
  The next illustration represents a sample of the random tree
  $\mathcal{T}_{t,x,{\rm Id}}$ started from $c = {\rm Id}$
  for a semilinear PDE of the form \eqref{eq:112}. 
  
\tikzstyle{level 1}=[level distance=3cm, sibling distance=1cm]
\tikzstyle{level 2}=[level distance=4cm, sibling distance=4cm]

\vspace{-0.5cm}

\begin{center}
\resizebox{0.85\textwidth}{!}{
\begin{tikzpicture}[yscale = 0.55,scale=1.3,grow=right, sloped,color=blue]
\node[rectangle,draw,black,text=black,thick]{$t$}
            child {
              node[rectangle,draw,black,text=black,thick] {$T_{(1)}$}
              child {
                node[rectangle,draw,black,text=black,thick] {$T_{(1,1)}$}
                child{
                node[rectangle,draw,black,text=black,thick]{$T_{(1,1,2)}$}
                    child{
                node[rectangle,draw,black,text=black,thick]{$T_{(1,1,2,2)}$}
                child{
                    node[rectangle,draw,black,text=black,thick]{$T$}
                    edge from parent
                    node[above]{$(1,1,2,2,2)$}
                    node[below]{$(f^{(3)})^*$}
                    }
                    child{
                    node[rectangle,draw,black,text=black,thick]{$T$}
                    edge from parent
                    node[above]{$(1,1,2,2,1)$}
                    node[below]{$(f^{(2)})^*$}
                    }
                    edge from parent
                node[above]{$(1,1,2,2)$}
                node[below]{$( f^{(2)})^2$}
                }
                    child{
                    node[rectangle,draw,black,text=black,thick]{$T$}
                    edge from parent
                    node[above]{$(1,1,2,1)$}
                    node[below]{$(f')^*$}
                    }
                edge from parent
                node[above]{$(1,1,2)$}
                node[below]{$(f')^*$}
                }
                child{
                node[above=10pt,rectangle,draw,black,text=black,thick,yshift=2cm]{$T_{(1,1,1)}$}
                    child{
                    node[rectangle,draw,black,text=black,thick]{$T$}
                    edge from parent
                    node[above]{$(1,1,1,3)$}
                    node[below,yshift=-0.3cm]{$(-(1/2) f^{(2)} )^*$}
                    }
                    child{
                    node[rectangle,draw,black,text=black,thick]{$T$}
                    edge from parent
                    node[above]{$(1,1,1,2)$}
                    node[below]{$\partial_x$}
                    }
                    child{
                    node[rectangle,draw,black,text=black,thick]{$T$}
                    edge from parent
                    node[above]{$(1,1,1,1)$}
                    node[below]{$\partial_x$}
                    }
edge from parent
                node[above]{$(1,1,1)$}
                node[below]{$f^*$}
                                }
                edge from parent
                node[above] {$(1,1)$}
                node[below]  {$f^*$}
            }
                edge from parent
                node[above] {$(1)$}
                node[below]  {${\rm Id}$}
            };
\end{tikzpicture}
}
\end{center}

\section{Probabilistic representation of PDE solutions} 
\label{s3}
 In this section, we derive a probabilistic representation formula
 for the solution of fully nonlinear PDEs of the form \eqref{eq:1}, 
 using a multiplicative functional $\mathcal{H}(\mathcal{T}_{t,x,c})$
 of the random tree $\mathcal{T}_{t,x,c}$. 
 For this, we will link the codes introduced in the previous section
 to the Duhamel formulation of the PDE \eqref{eq:1} 
 by deriving a system of equations satisfied by $\E [\mathcal{H}(\mathcal{T}_{t,x,c}) ]$,
 $c\in \eufrak{C}$. 
 We let 
 $$
 \widebar{F}(t) := \int_t^\infty \rho(u) du = \mathbb{P}(\tau > t), \qquad t \in \real_+, 
 $$
 denote the tail distribution function of $\tau^{i,j}$, $i,j\geq 1$.  
\begin{definition}
 We define the functional $\mathcal{H}(\mathcal{T}_{t,x,c})$ 
 of the random coding tree $\mathcal{T}_{t,x,c}$
 started at time $t\in [0,T]$, location $x\in \real$
 and code $c \in \eufrak{C}$ as 
 $$
 \mathcal{H}(\mathcal{T}_{t,x,c}) := \prod_{\widebar{k} \in \mathcal{K}^{\circ}} \frac{1}{q_{c_{\widebar{k}}}(I_{c_{\widebar{k}}})\rho(\tau_{\widebar{k}})} \prod_{\widebar{k} \in \mathcal{K}^{\partial}} \frac{c_{\widebar{k}}(u)\big(T,X^{\bar{k}}_{T,x} \big)}{\widebar{F}(T-T_{\widebar{k}-})}.
   $$
\end{definition}{}
Note that for $c \in \eufrak{C}$ of the form $c= \partial_x^k$ we have 
$$
c (u)(T, x ) = \partial_x^k \phi ( x), 
$$
and for $c\in \eufrak{C}$ of the form $c = g^*$ with
$g = a \partial_{z_0}^{\lambda_0} \cdots \partial_{z_n}^{\lambda_n} f$, we have 
$$
c (u)(T, x ) = g\big(\phi (T,x),\partial_x\phi (T,x), \ldots , \partial_x^n\phi (T,x)\big), \quad ( \lambda_0,\ldots , \lambda_n ) \in \inte^{n+1}.
$$ 
The next result gives the probabilistic representation of solutions
of \eqref{eq:1} as an expectation over random coding trees.
\begin{theorem}
   \label{t1}
    Under Assumption~(\hyperlink{BGJhyp}{$A$}), let $T>0$ such that 
   \begin{equation}
\nonumber 
   \E \big[ \big| \mathcal{H}(\mathcal{T}_{t,x,c}) \big| \big]
 < \infty,
 \quad (t,x) \in [0,T] \times \real,
  \ c\in \eufrak{C}. 
\end{equation} 
Then, for any $c\in \eufrak{C}$ the function 
 \begin{equation}
 \nonumber 
   u_c(t,x): = \E \big[ \mathcal{H}(\mathcal{T}_{t,x,c}) \big], \qquad
   (t,x)\in [0,T]\times \real, 
 \end{equation}
 is a solution of the system 
\begin{equation} 
\label{s1-1}
 u_c(t,x) = \int_{-\infty}^\infty \varphi (T-t,y-x) c(u)(T,y) dy 
+ 
\sum_{Z \in \mathcal{M}(c)}
\int_t^T \int_{-\infty}^\infty \varphi (s-t,y-x)
\prod_{z \in Z}  u_z (s,y) dy ds, 
\end{equation}
 with $u_c (T,x) = c(u)(T,x)$, $(t,x)\in [0,T] \times \real$.
 Moreover, if the solution $(u_c)_{c\in \eufrak{C}}$
of \eqref{s1-1} is unique, then we have 
$$
c(u) (t,x) = u_c(t,x) = \E \big[ \mathcal{H}(\mathcal{T}_{t,x,c}) \big], \qquad
   (t,x) \in [0,T] \times \real. 
$$ 
  In particular, taking $c={\rm Id}$, we have the probabilistic representation 
 \begin{equation}
 \label{fjhkldsf} 
   u(t,x) = \E \big[ \mathcal{H}(\mathcal{T}_{t,x , {\rm Id}}) \big],
    \qquad
   (t,x) \in [0,T] \times \real.   
\end{equation}
\end{theorem}
\begin{Proof}
 For $c\in \eufrak{C}$, we let 
$$
   u_c(t,x) := \E \big[ \mathcal{H}(\mathcal{T}_{t,x,c})\big],
   \qquad
   (t,x) \in [0,T] \times \real.  
$$
 By conditioning on the first branching time $T_{\widebar{1}}$, 
 the first particle bearing the code ${\rm Id}$
branches at time $T_{\widebar{1}}$
into a new particle bearing the code $f^*$
as $\mathcal{M} ( {\rm Id} ) = \{ f^* \}$, hence 
\begin{eqnarray*}
  \lefteqn{
    u_{\rm Id}(t,x) =  \E \big[ \mathcal{H}(\mathcal{T}_{t,x,{\rm Id}})\mathbbm{1}_{\{ T_{\widebar{1}}> T \}} + \mathcal{H}(\mathcal{T}_{t,x,{\rm Id}}) \mathbbm{1}_{\{ T_{\widebar{1}}\leq T \}} \big]
  }
  \\
&= & \E \Bigg[ \frac{\phi \big( X^{\widebar{1}}_{T,x}\big)}{\widebar{F}(T-t)} \mathbbm{1}_{\{ T_{\widebar{1}}>T \}} \Bigg]
+  \E \Bigg[
  \frac{u_{f^*}\big( T_{\widebar{1}}, X^{\widebar{1}}_{T_{\widebar{1}},x}\big)}{\rho(T_{\widebar{1}} - t )}\mathbbm{1}_{\{ T_{\widebar{1}} \leq T \}}  \Bigg]\\
&= & \frac{\mathbb{P}(T_{\widebar{1}}>T)}{ \widebar{F}(T-t)}
\int_{-\infty}^\infty \varphi (T-t,y-x) \phi (y) dy
+ \int_t^T
\int_{-\infty}^\infty \varphi (s-t,y-x) u_{f^*}(s,y) dy ds
\\
&= & \int_{-\infty}^\infty \varphi (T-t,y-x) \phi (y) dy
+ \int_t^T
\int_{-\infty}^\infty \varphi (s-t,y-x) u_{f^*}(s,y) dy ds,
\quad
(t,x) \in [0,T] \times \real.  
\end{eqnarray*} 
Similarly, starting from any code $c \in \eufrak{C}$ different from ${\rm Id}$,
we draw a sample of $I_c$ uniformly in $\mathcal{M}(c)$.
As each code in the tuple $I_c$ yields a new branch 
at time $T_{\widebar{1}}$, we obtain 
\begin{align}
\nonumber
&
\! \! \! \! \! \! 
 u_c (t,x) =  \E \big[ \mathcal{H}(\mathcal{T}_{t,x,c})\mathbbm{1}_{\{ T_{\widebar{1}}> T \}} + \mathcal{H}(\mathcal{T}_{t,x,c}) \mathbbm{1}_{\{ T_{\widebar{1}}\leq T \}} \big]
\\
\nonumber
&= \E \left[ \frac{{c(u)\big(T, X^{\widebar{1}}_{T,x}\big) }}{\widebar{F}(T-t)} \mathbbm{1}_{\{ T_{\widebar{1}}> T \}}
+
  \mathbbm{1}_{\{ T_{\widebar{1}} \leq T \}}
  \sum_{Z\in {\cal M}(c)}
  \mathbbm{1}_{\{ I_c = Z\}} 
\frac{\prod_{z\in Z} u_z \big( T_{\widebar{1}}, X^{\widebar{1}}_{T_{\widebar{1}},x}\big)}{ q_c( Z ) \rho(T_{\widebar{1}} - t )}
\right]
\\
\nonumber 
  & =
  \int_{-\infty}^\infty \varphi (T-t,y-x) c(u)(T,y) dy 
+ 
\sum_{Z \in \mathcal{M}(c)}
\int_t^T \int_{-\infty}^\infty \varphi (s-t,y-x)
\prod_{z \in Z}  u_z (s,y) dy ds, 
\end{align}
which yields the system of equations \eqref{s1-1}.
 We conclude by noting that from Lemma~\ref{l1}, the family of functions  
 $(c (u))_{c\in \eufrak{C} }$ is the solution of the system \eqref{s1-1}, 
 hence we have 
 $(c (u))_{c\in \eufrak{C} } = (u_c)_{c\in \eufrak{C} }$, and 
$$
\E \big[ \mathcal{H}(\mathcal{T}_{t,x,c}) \big]
=
u_c (t,x) = 
c(u)(t,x), \qquad
(t,x) \in [0,T] \times \real,
\quad c \in \eufrak{C}.  
$$
\end{Proof} 
 In the case of a semilinear PDE of the form \eqref{eq:112} with $n=0$, 
 the system
\eqref{s1-1}
 reads 
 \begin{equation}
 \label{eq:2-2}
\begin{cases}
  \displaystyle
  \partial_t u_{{\rm Id}}(t,x) + \frac{1}{2}\partial_x^2 u_{{\rm Id}}(t,x) + u_{f^*}(t,x) = 0
  \medskip
  \\
\displaystyle
\partial_t u_{( f^{(k)} )^*}(t,x) + \frac{1}{2}\partial_x^2 u_{( f^{(k+2)} )^*}(t,x) + u_{( f^{(k)} )^*}(t,x) u_{( f^{(k+1)} )^*}(t,x) - \frac{1}{2}( u_{\partial_x} (t,x) )^2 u_{( f^{(k+2)} )^*}(t,x) = 0, 
\medskip 
  \\
\displaystyle
\partial_t u_{\partial_x}(t,x) + \frac{1}{2}\partial_x^2 u_{\partial_x}(t,x) + u_{\partial_x}(t,x) u_{( f')^*} (t,x) = 0
\medskip
  \\
\displaystyle
u_c(T,x) = c (u) (T,x), \quad c\in \eufrak{C}, 
\end{cases}
\end{equation} 
 $k \geq 0$. 
 For example, in case $n=0$ and $f(z_0) = \re^{z_0}$,
 the system \eqref{eq:2-2} simplifies to 
 $$
\begin{cases}
  \displaystyle
  \partial_t u_{{\rm Id}}(t,x) + \frac{1}{2}\partial_x^2 u_{{\rm Id}}(t,x) + u_{f^*}(t,x) = 0
  \medskip
  \\
\displaystyle
\partial_t u_{f^*} (t,x) + \frac{1}{2}\partial_x^2 u_{f^*} (t,x) + ( u_{f^*} (t,x) )^2 - \frac{1}{2}( u_{\partial_x} (t,x) )^2 u_{f^*} (t,x) = 0, 
\medskip 
  \\
\displaystyle
\partial_t u_{\partial_x}(t,x) + \frac{1}{2}\partial_x^2 u_{\partial_x}(t,x) + u_{\partial_x}(t,x) u_{f^*} (t,x) = 0
\medskip
  \\
\displaystyle
u_c(T,x) = c (u) (T,x), \quad c\in \eufrak{C}. 
\end{cases}
$$
The next result provides sufficient conditions for
the uniform boundedness of the random functional $\mathcal{H}(\mathcal{T}_{t,x,c})$,
therefore ensuring the integrability needed for the
validity of the probabilistic representation \eqref{fjhkldsf} in Theorem~\ref{t1}. 
 \begin{prop} 
 \label{t3.3}
 Under Assumption~(\hyperlink{BGJhyp}{$A$}),
 suppose in addition that 
 $$
  \Vert \partial_x^k \phi \Vert_{L^\infty (\real)} \leq K, \
  \big\Vert \partial_{z_0}^{\lambda_0} \cdots \partial_{z_n}^{\lambda_n} f\big(\phi , \partial_x \phi , \ldots , \partial_x^n \phi \big) \big\Vert_{L^\infty (\real)} \leq K,
  \quad k\geq 0 , \ \lambda_0,\ldots , \lambda_n \geq 0, 
$$ 
 for some $K \in (0,1)$, and that the probability density function $\rho (t)$ is 
decreasing and satisfies the conditions 
 $$
\rho(T) \geq \frac{1}{\min_{c \in \eufrak{C}} q_c(I_c)}
\quad
 \mbox{and} \quad 
K \leq \widebar{F}(T) 
$$
for some $T>0$.
Then we have $| \mathcal{H}(\mathcal{T}_{t,x,c} )|\leq 1$, $a.s.$,
$(t,x,c)\in [0,T]\times \real\times \eufrak{C}$. 
\end{prop} 
\begin{Proof}
 Since 
 $\partial_x^k (u)(T, x )$ and $g^* (u)(T, x )$ respectively denote
 the functions $\partial_x^k \phi ( x)$ and 
 $g\big(\phi (T,x),\partial_x\phi (T,x), \ldots , \partial_x^n\phi (T,x)\big)$,
 we have the bound  
\begin{eqnarray*} 
  | \mathcal{H}(\mathcal{T}_{t,x,c} )|
  & \leq & 
  \prod_{\widebar{k} \in \mathcal{K}^{\circ}} \frac{1}{q_{c_{\widebar{k}}}(I_{c_{\widebar{k}}})\rho(\tau_{\widebar{k}})} \prod_{\widebar{k} \in \mathcal{K}^{\partial}} \frac{K}{\widebar{F}(T-T_{\widebar{k}-})}
  \\
   & \leq & 
    \prod_{\widebar{k} \in \mathcal{K}^{\circ}} \frac{1}{\rho(T)
      \min_{c \in \eufrak{C}} q_c(I_c) 
    } \prod_{\widebar{k} \in \mathcal{K}^{\partial}} \frac{K}{\widebar{F}(T)}
    \\
    & \leq & 1, \quad a.s.. 
\end{eqnarray*} 
\end{Proof}
\section{Numerical examples}
\label{s5}
 In this section we provide numerical confirmations of the validity
 of our algorithm on examples of nonlinear and fully nonlinear PDEs, by benchmarking
 its output to closed-form solutions. 
 Numerical computations are done in Mathematica using
 the codes provided in appendix, which apply to the general
 fully nonlinear case with higher order derivatives, 
 with $\rho$ the standard exponential probability density.
 Semilinear examples are treated using a simplified code that
 does not use gradient nonlinearities. 
\subsubsection*{Semilinear examples} 
\begin{description} 
\item{Example~1-$a)$.}
 Consider the Allen-Cahn (or Ginzburg-Landau) equation
\begin{equation}
\label{gl0} 
\partial_t u(t,x) + \Delta_x u(t,x) + u(t,x) - u^3(t,x) = 0, 
\quad
(t,x) \in [0,T]\times \real^d, 
\end{equation} 
as in \cite{han2018solving} or \S~4.2 of \cite{han2018solvingarxiv},
with terminal condition $\phi (x) = 1/(2 + 2 \Vert x\Vert^2 / 5)$, 
where $\Vert x\Vert := \sqrt{x_1^2+\cdots + x_d^2}$,
$x=(x_1,\ldots , x_d) \in \real^d$. 
In Table~\ref{fjkldsf1} we compare our results to the ones obtained
in Table~1 of \cite{han2018solvingarxiv}
for the estimation of $u(0,x)$ at $x=(0,\ldots , 0)$ for $T=0.3$ 
with 4000 iterations in dimension $d=100$. 
We note that such results can be recovered by the multilevel Picard method,
see also \S~3 of \cite{hutzenthaler-mlp4}. 
 
\begin{table}[H]
  \centering
\footnotesize 
    \begin{tabular}{|l|c|c|}
   \hline
   & \cite{han2018solvingarxiv} & Coding trees 
   \\
   \hline
 Mean & 0.0528 & 0.052754
 \\
 \hline
 Standard deviation & 0.0002 & 0.000364
  \\
  \hline
  Mean of rel. $L^1$ error & 0.0030 & 0.005916
  \\
  \hline
  SD of rel. $L^1$ error & 0.0022 & 0.003661
 \\
 \hline
\end{tabular}
    \caption{BSDE \cite{han2018solvingarxiv} {\em vs} coding tree method for \eqref{gl0}
 with $T=0.3$ and $d=100$. 
}  
\label{fjkldsf1}
\end{table}

\vspace{-0.55cm}

Next, we consider the Allen-Cahn equation
\begin{equation}
\label{gl} 
\partial_t u(t,x) + \frac{1}{2} \Delta_x u(t,x) + u(t,x) - u^3(t,x) = 0, 
\quad
(t,x) \in [0,T]\times \real^d, 
\end{equation} 
 which admits the traveling wave solution
\begin{equation}
  \label{fjklds}
  u(t,x) = -\frac{1}{2} - \frac{1}{2}
\tanh \left( \frac{3}{4} (T-t) - \sum_{i=1}^d \frac{x_i}{2 \sqrt{d}} \right),
\quad
(t,x) \in [0,T]\times \real^d. 
\end{equation}
In Figure~\ref{fig1} we compare the closed-form solution \eqref{fjklds} of
\eqref{gl} to its estimation by the coding tree method. 
 
\begin{figure}[H]
\centering
\hskip0.2cm
\begin{subfigure}{.49\textwidth}
  \hskip0.3cm
  \includegraphics[width=\textwidth]{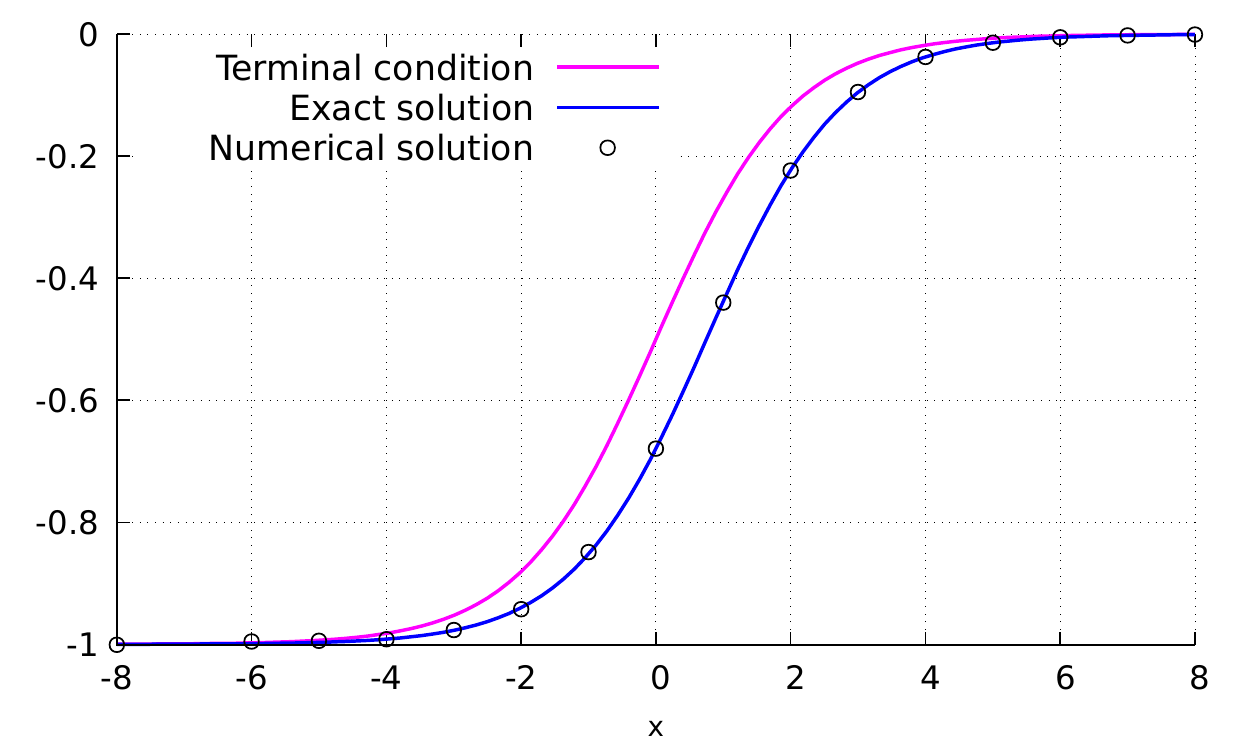}
\vskip-0.1cm
\caption{Dimension $d=1$ with $T=0.5$.}
\end{subfigure}
\begin{subfigure}{.49\textwidth}
\centering
\includegraphics[width=\textwidth]{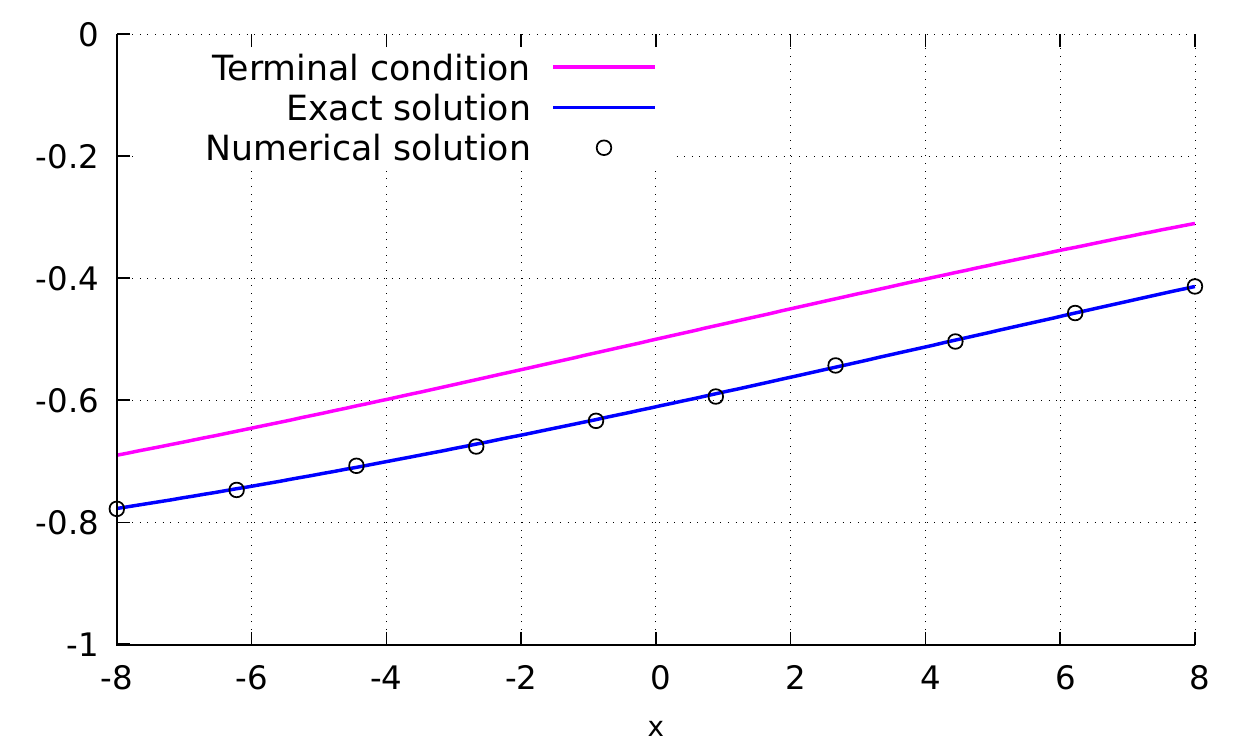}
\vskip-0.1cm
\caption{Dimension $d=100$ with $T=0.3$.}
\end{subfigure}
\caption{Numerical solution $u(0,x)$ of \eqref{gl} with $100,000$ Monte Carlo samples.} 
\label{fig1}
\end{figure}

\vspace{-0.5cm}

In Table~\ref{t3} we compare the stability of the coding tree method 
to that of the branching diffusion of \cite{laborderespa} according
to Table~8 in  \cite{hutzenthaler-mlp0} in dimension $d=1$
 with $T=1$ and the explicit solution 
 \begin{equation}
   \label{fjklf} 
 u(t,x) = \frac{1}{\sqrt{1-(1-(\phi (0))^{-2}) \re^{-2(T-t)}}}, \qquad t\in [0,T]. 
\end{equation} 
 It turns out that our coding tree algorithm remains stable for higher values of
 $\phi (0)$ in this example. 
 
\begin{table}[H]
    \centering
\footnotesize 
    \begin{tabular}{|c|c|c|c}
   \hline
 $\phi (0)$   & Exact value $u(0,0)$ & Coding trees & \multicolumn{1}{c|}{\cite{laborderespa}}  
  \\
  \hline
0.1 & 0.263540 & 0.247403 & \multicolumn{1}{c|}{0.271007} 
  \\
  \hline
0.2 & 0.485183 & 0.472720 & \multicolumn{1}{c|}{0.499103}  
  \\
  \hline
0.3 & 0.649791 & 0.723543 & \multicolumn{1}{c|}{0.848879} 
  \\
  \hline
0.4 & 0.764605 & 0.866281 & \multicolumn{1}{c|}{3.495457} 
  \\
  \hline
0.5 & 0.843347 & 0.932852 & \multicolumn{1}{c|}{21.68436} 
  \\
  \hline
0.6 & 0.897811 & 0.968213 & \multicolumn{1}{c|}{136.6667} 
  \\
  \hline
0.7 & 0.936233 & 1.005440 & \multicolumn{1}{c|}{7321.326} 
  \\
  \hline
0.8 & 0.963981 & 0.950816 & 
  \\
  \cline{1-3}
0.9 & 0.984496  & 0.944715 & 
  \\
  \cline{1-3}
1.0 & 1.0  & 1.000164 & 
\\    \cline{1-3} 1.1 & 1.011955 & 1.182766 & 
 \\    \cline{1-3} 1.2 & 1.021340  & 1.576551 & 
  \\
  \cline{1-3}
1.5 & 1.039856 & 5.182978 & 
  \\
  \cline{1-3}
2.0 & 1.054973 & 30.006351 & 
\\
  \cline{1-3}
\end{tabular}
\caption{Branching diffusion \cite{laborderespa} {\em vs} coding trees for the Allen-Cahn equation \eqref{gl}.} 
\label{t3}
\end{table}

\vspace{-0.6cm} 
 
The stability over time of the coding tree algorithm applied to the
Allen-Cahn equation \eqref{gl} with solution \eqref{fjklf}
is illustrated in Figure~\ref{f0-exp}. 

\begin{figure}[H]
\centering
\hskip0.2cm
\begin{subfigure}{.49\textwidth}
  \hskip0.3cm
  \includegraphics[width=\textwidth]{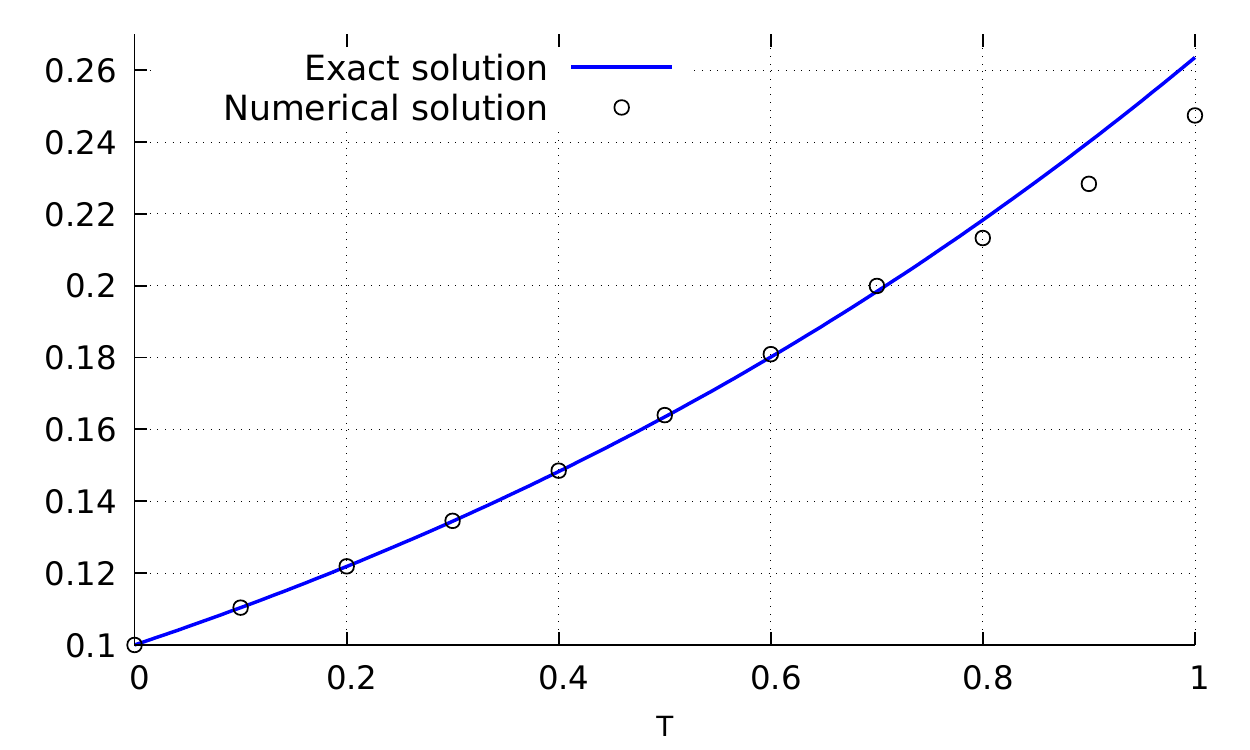}
\vskip-0.1cm
\caption{$\phi (0)=0.1$.}
\end{subfigure}
\begin{subfigure}{.49\textwidth}
\centering
\includegraphics[width=\textwidth]{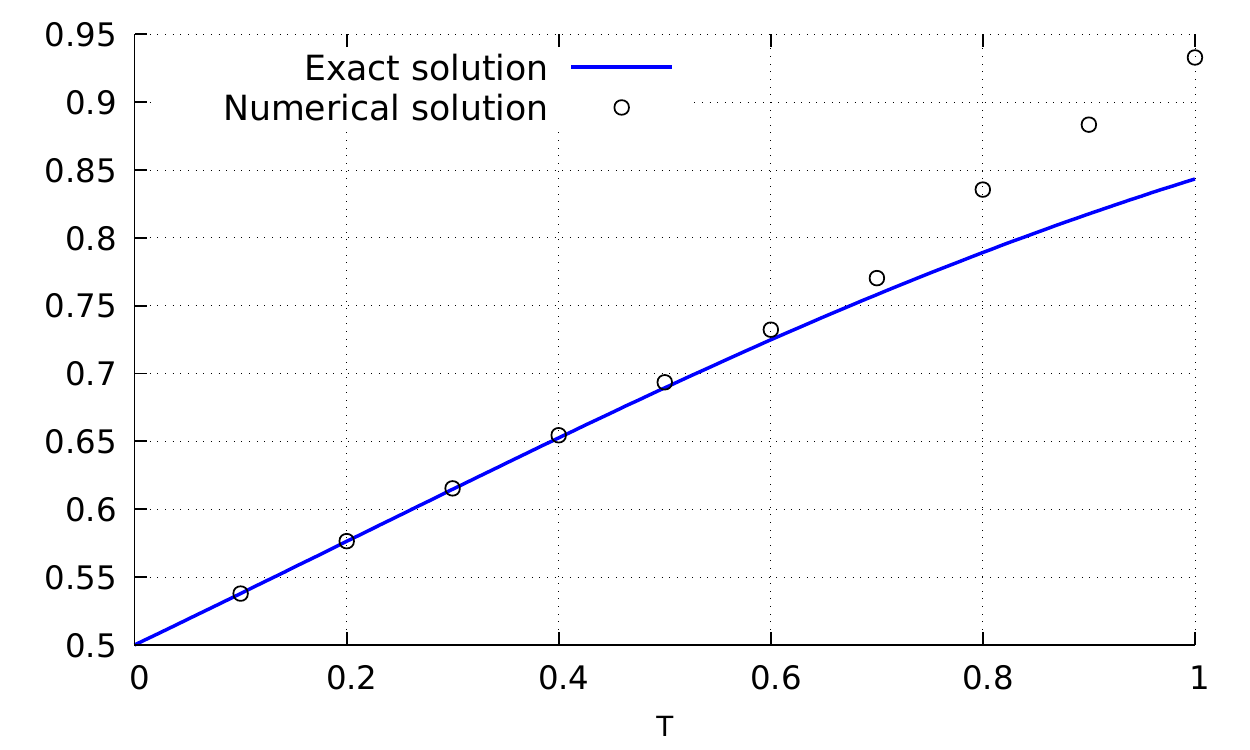}
\vskip-0.1cm
\caption{$\phi (0)=0.5$.}
\end{subfigure}
\caption{Numerical solution $u(0,0)$ of \eqref{gl} by the coding tree method.}
\label{f0-exp}
\end{figure}

\vskip-0.3cm

In Figure~\ref{ac}, the stability over time of the coding tree algorithm is compared to that of the BSDE method \cite{han2018solvingarxiv}
 for the Allen-Cahn equation \eqref{gl} with solution \eqref{fjklds}, 
 using the 
 BSDE solver available at \url{https://github.com/frankhan91/DeepBSDE}. 

\begin{figure}[H]
\centering
\hskip0.2cm
\begin{subfigure}{.49\textwidth}
  \hskip0.3cm
\includegraphics[width=\textwidth]{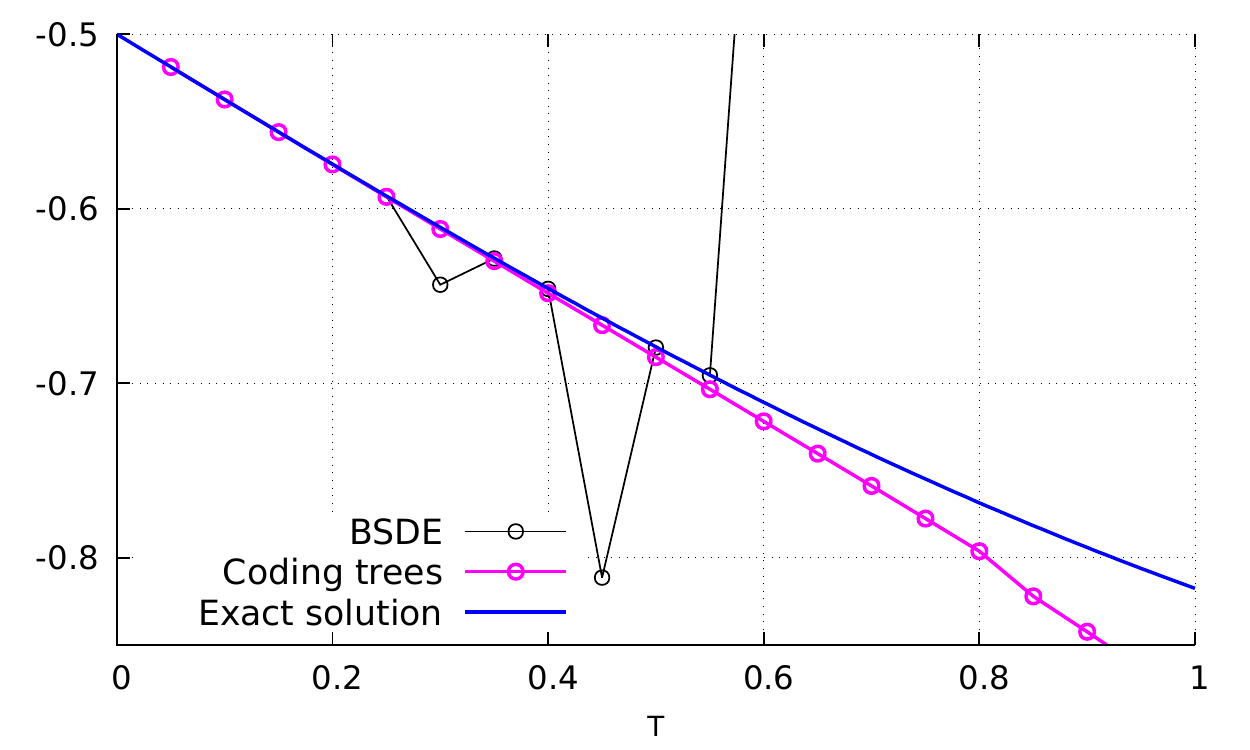}
\vskip-0.1cm
\caption{Dimension $d=1$.}
\end{subfigure}
\begin{subfigure}{.49\textwidth}
\centering
\includegraphics[width=\textwidth]{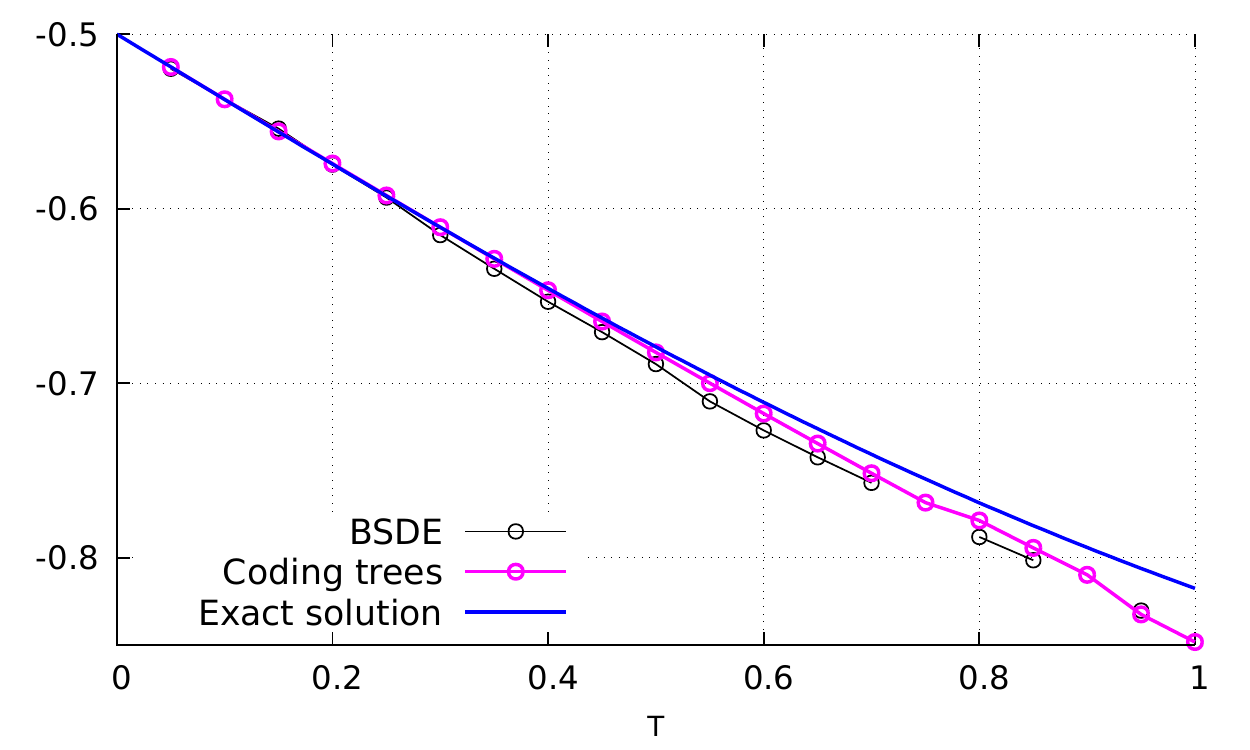}
\vskip-0.1cm
\caption{Dimension $d=100$.}
\end{subfigure}
    \caption{      BSDE method \cite{han2018solvingarxiv} {\em vs} coding trees for the Allen-Cahn equation \eqref{gl}.} 
\label{ac} 
\end{figure}

\vskip-0.5cm

The data of Figure~\ref{ac} with $d=100$ is presented in Table~\ref{t3-0},
 where it turns out that our coding tree algorithm remains stable for higher values of
 $T$ and that the BSDE method is less stable in low dimension in this example.
 We also note that this BSDE solver requires the input of
 an initial guess interval y\_init\_range for the algorithm to run,
 which is not the case in branching type methods. 

\begin{table}[H]
    \centering
    \footnotesize 
    \begin{tabular}{|c|c|c|c|c|c|}
  \cline{3-6}
 \multicolumn{1}{c}{}  &  \multicolumn{1}{c}{}   & \multicolumn{2}{|c|}{$d=1$} & \multicolumn{2}{|c|}{$d=100$}
       \\
   \hline
 $T$ & Exact value & Coding trees & BSDE & Coding trees & BSDE 
  \\
  \hline
   0.10 & -0.537430  & -0.537451 &   -0.537470   & -0.537318  & -0.537169
  \\
  \hline
   0.20 & -0.574443  & -0.574662  &  -0.574581  & -0.574111  & -0.574704
  \\
  \hline
   0.30  & -0.610639 & -0.611628  &  -0.643579  & -0.610616  & -0.615149
  \\
  \hline
   0.40 & -0.645656  & -0.648401  &  -0.645983   & -0.646737  & -0.653324
  \\
  \hline
   0.50 & -0.679179  & -0.685096  &  -0.679532   & -0.682355  & -0.688992
  \\
  \hline
  0.60 & -0.710950  & -0.721866  &  -0.270275    & -0.717349  & -0.727000
  \\
  \hline
  0.70 & -0.740775  & -0.758880     &     NaN    & -0.751602  & -0.757082
  \\
  \hline
  0.80 & -0.768525  & -0.796317      &    NaN    & -0.778700  & -0.788100
  \\
  \hline
  0.90 & -0.794130 & -0.842476      &    NaN  & -0.809848  &       NaN 
  \\
  \hline
  1.00 & -0.817574  & -0.884167      &    NaN    & -0.848294  &       NaN
  \\
  \hline
  1.20 & -0.858149  & -0.945912      &    NaN   & -0.900085  & -0.530944
  \\
  \hline
  1.40 & -0.890903  & -1.023799       &   NaN    & -0.965941  & -0.473223
  \\
  \hline
  1.60 & -0.916827  & -1.104303       &   NaN    & -1.021344  &       NaN
  \\
  \hline
  1.80 & -0.937027  & -1.155304       &   NaN     & -1.089572  &       NaN
  \\
  \hline
  2.00  & -0.952574 & -1.193363       &   NaN     & -1.127749  &       NaN
  \\
  \hline
    \end{tabular}
    \caption{      BSDE method \cite{han2018solvingarxiv} {\em vs} coding trees for the Allen-Cahn equation \eqref{gl}.} 
\label{t3-0}
\end{table}

\vspace{-0.6cm} 

\noindent
\item{Example~1-$b)$.}
  As the Allen-Cahn Example~1-$a)$ only involves polynomial nonlinearities,
  it can be treated by the branching diffusion method, see 
  \cite{henry-labordere2012}, \cite{labordere}. 
  On the other hand, the following example, which
  makes use of a functional nonlinearity, cannot be treated by such a method.
  Consider the equation 
\begin{equation}
\label{nonl0} 
\partial_t u(t,x) +
\frac{\alpha}{d} \sum\limits_{i=1}^d \partial_{x_i} u(t,x)
+ \frac{1}{2} \Delta_x u(t,x)
 + \re^{-u(t,x)} ( 1 - 2 \re^{-u(t,x)} ) d = 0, 
\end{equation} 
which admits the traveling wave solution
\begin{equation}
  \label{fjklds34} 
u(t,x) = \log \left( 1 + \left( \alpha (T-t) + \sum_{i=1}^d x_i \right)^2 \right), 
\qquad (t,x) \in [0,T]\times \real^d.
\end{equation} 
 In Figure~\ref{f0} we take $T=0.05$, $\alpha =10$ and $100,000$ Monte Carlo samples.
 
\begin{figure}[H]
\centering
\hskip0.2cm
\begin{subfigure}{.49\textwidth}
  \hskip0.3cm
  \includegraphics[width=\textwidth]{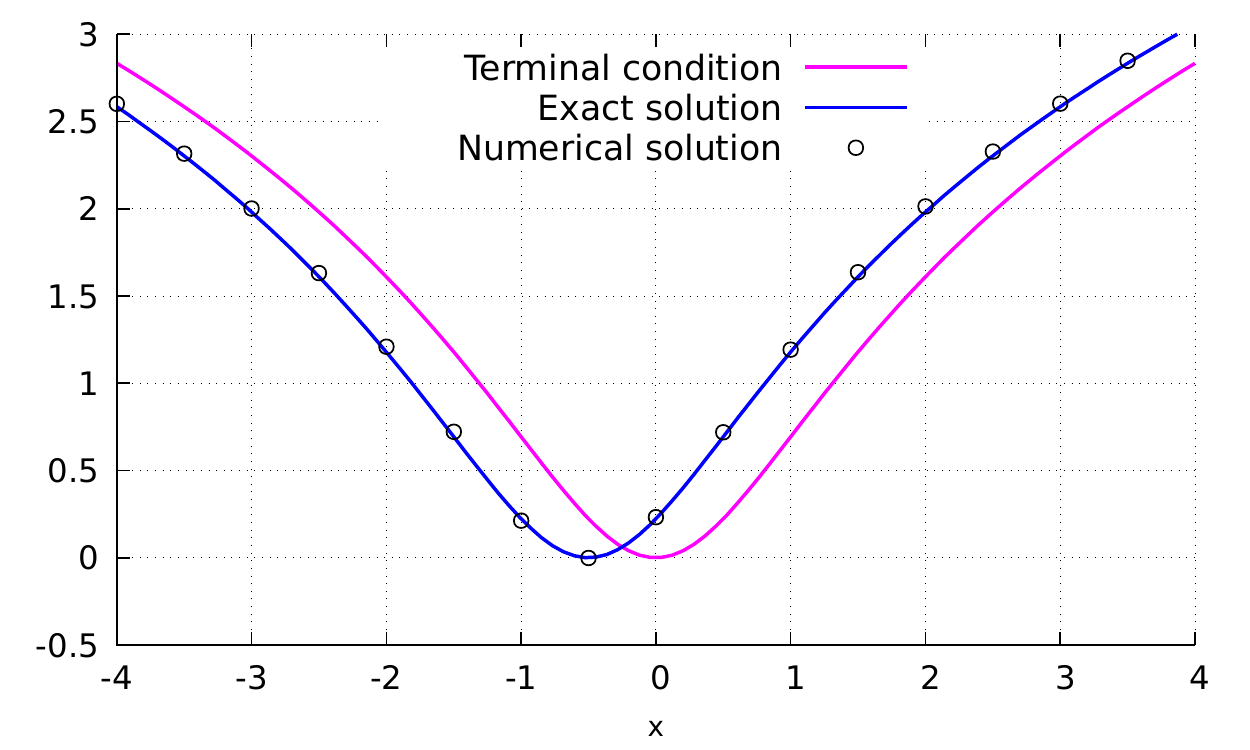}
\vskip-0.1cm
\caption{Dimension $d=1$.}
\end{subfigure}
\begin{subfigure}{.49\textwidth}
\centering
\includegraphics[width=\textwidth]{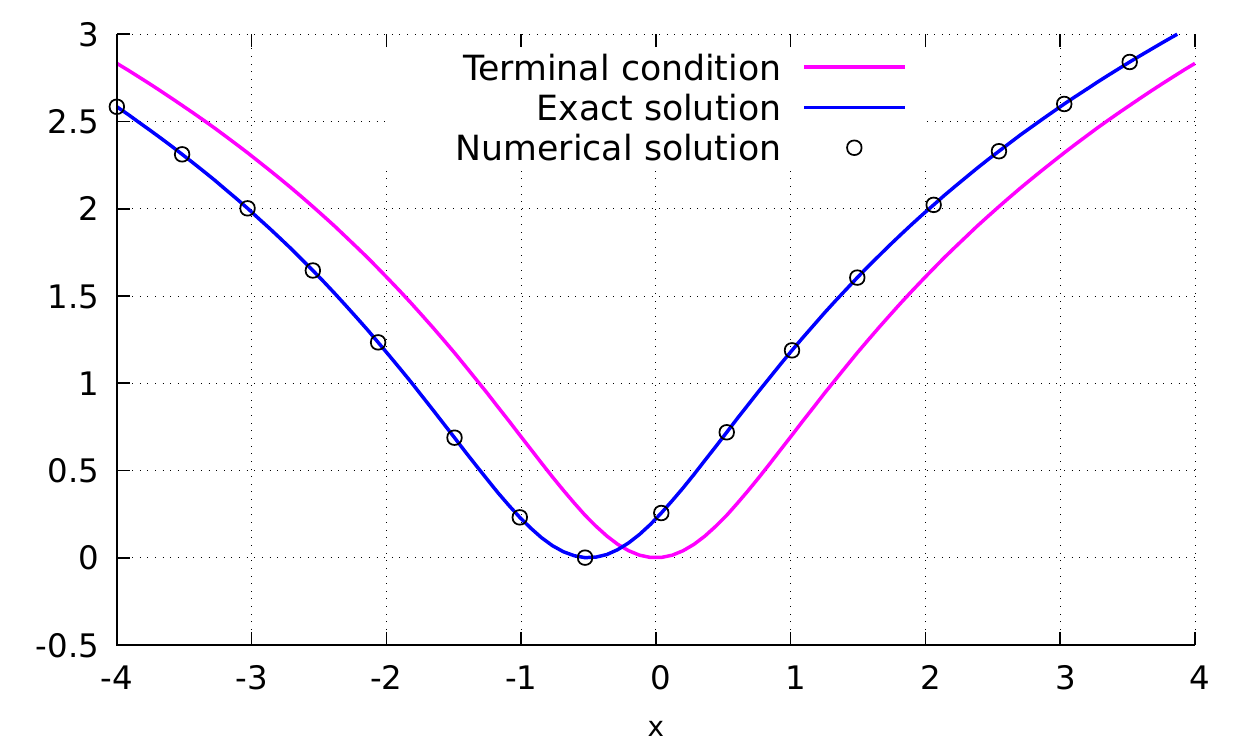}
\vskip-0.1cm
\caption{Dimension $d=10$.}
\end{subfigure}
\caption{Numerical solution $u(0,x)$ of \eqref{nonl0}.}
\label{f0}
\end{figure}

\vskip-0.3cm
 
 In Figure~\ref{f0-mlp} we take $T=0.05$, $\alpha =10$ and use the
 multilevel Picard code available at 
 \url{https://github.com/seb-becker/mlp} with $8$ iterations, 
 by including the following definition.

 \smallskip
 
 \begin{lstlisting}[language=C 
   ]
#define N_MAX 8
#define dimension 10
#define alphacoeff 10
#ifdef EXP_NONLIN
#define eq_name "Exponential nonlinearity"
#define rdim 1
#define TIME 0.05 
#define initial_value ArrayXd::Zero(d[j], 1)
#define g(x) ArrayXd tmp = ArrayXd::Zero(1, 1); tmp(0) = log ( 1 + pow(x.sum()+alphacoeff*TIME,2) )
#define X_sde(s, t, x, w) x + sqrt(1. * (t - s)) * w
#define fn(y) ArrayXd ret = ArrayXd::Zero(1, 1); double phi_r = std::min(4., std::max(-4., y(0))); ret(0) = exp(-phi_r)*(1-2*exp(-phi_r))*dimension
#endif
 \end{lstlisting}

 \vspace{-0.6cm} 
 
 We note that the performance of the multilevel Picard method is
 dimension-dependent in this example.

\begin{figure}[H]
\centering
\hskip0.2cm
\begin{subfigure}{.49\textwidth}
  \hskip0.3cm
  \includegraphics[width=\textwidth]{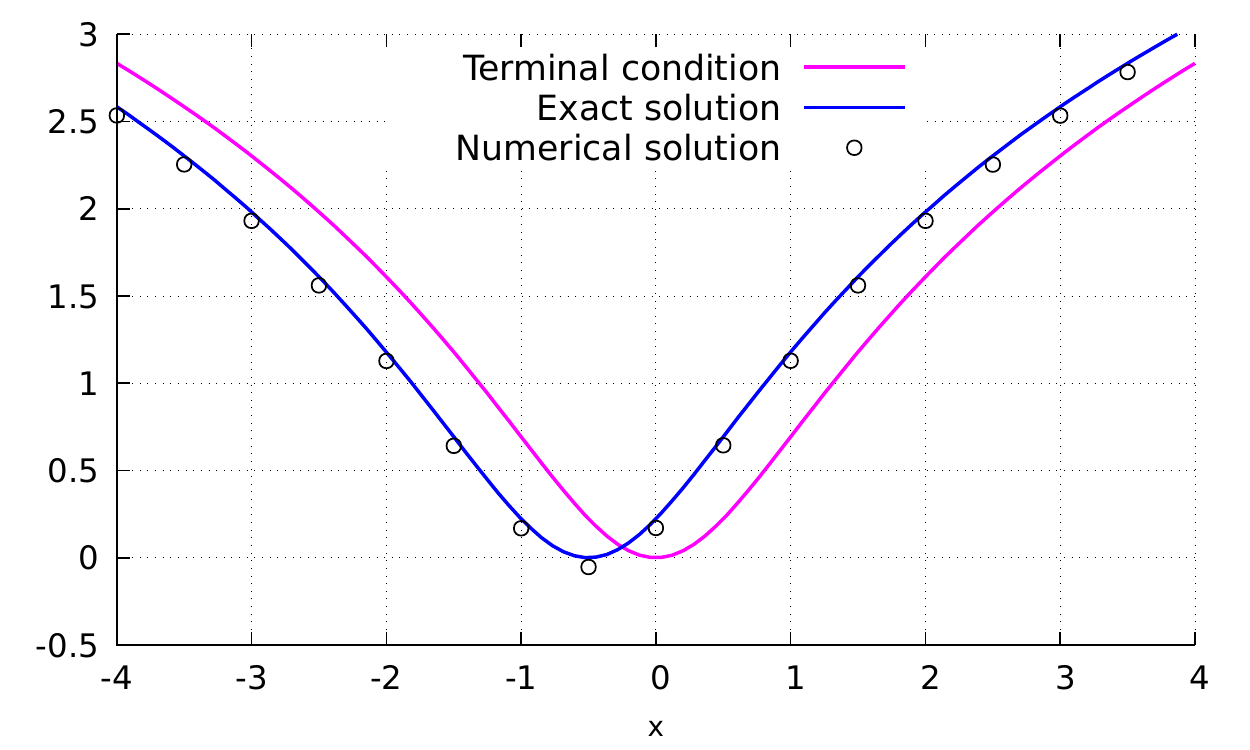}
\vskip-0.1cm
\caption{Dimension $d=1$.}
\end{subfigure}
\begin{subfigure}{.49\textwidth}
\centering
\includegraphics[width=\textwidth]{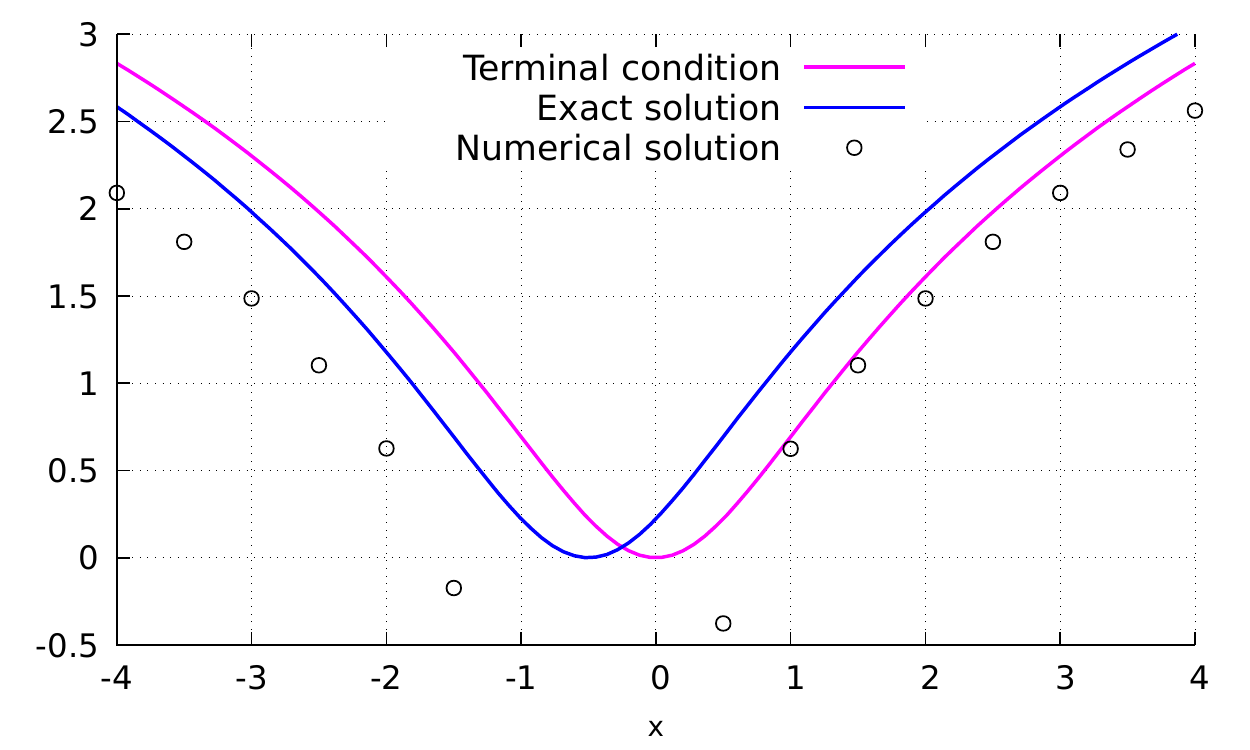}
\vskip-0.1cm
\caption{Dimension $d=10$.}
\end{subfigure}
\caption{Numerical solution $u(0,x)$ of \eqref{nonl0} by the multilevel Picard method.}
\label{f0-mlp}
\end{figure}

\vspace{-0.6cm}

\end{description}

\subsubsection*{Quasilinear examples}
\begin{description}
\item{Example~2-$a)$.}
 Consider the Dym equation
 \begin{equation}
   \label{dym} 
   \partial_t u(t,x) + \frac{1}{d} u^3(t,x)
   \sum\limits_{i=1}^d \partial_{x_i}^3 u(t,x)
   = 0, 
\end{equation} 
which admits the traveling wave solution
$$
u(t,x) =\left(3\alpha \left( 4 \alpha^2 (T - t) + \sum_{i=1}^d x_i  \right)\right)^{2/3},
\qquad
(t,x) \in [0,T]\times \real^d.
$$ 
with $\alpha > 0$. In Figure~\ref{f3} we take $T=0.01$, $\alpha=2$
and $100,000$ Monte Carlo samples.
 
\begin{figure}[H]
\centering
\hskip0.2cm
\begin{subfigure}{.49\textwidth}
  \hskip0.3cm
  \includegraphics[width=\textwidth]{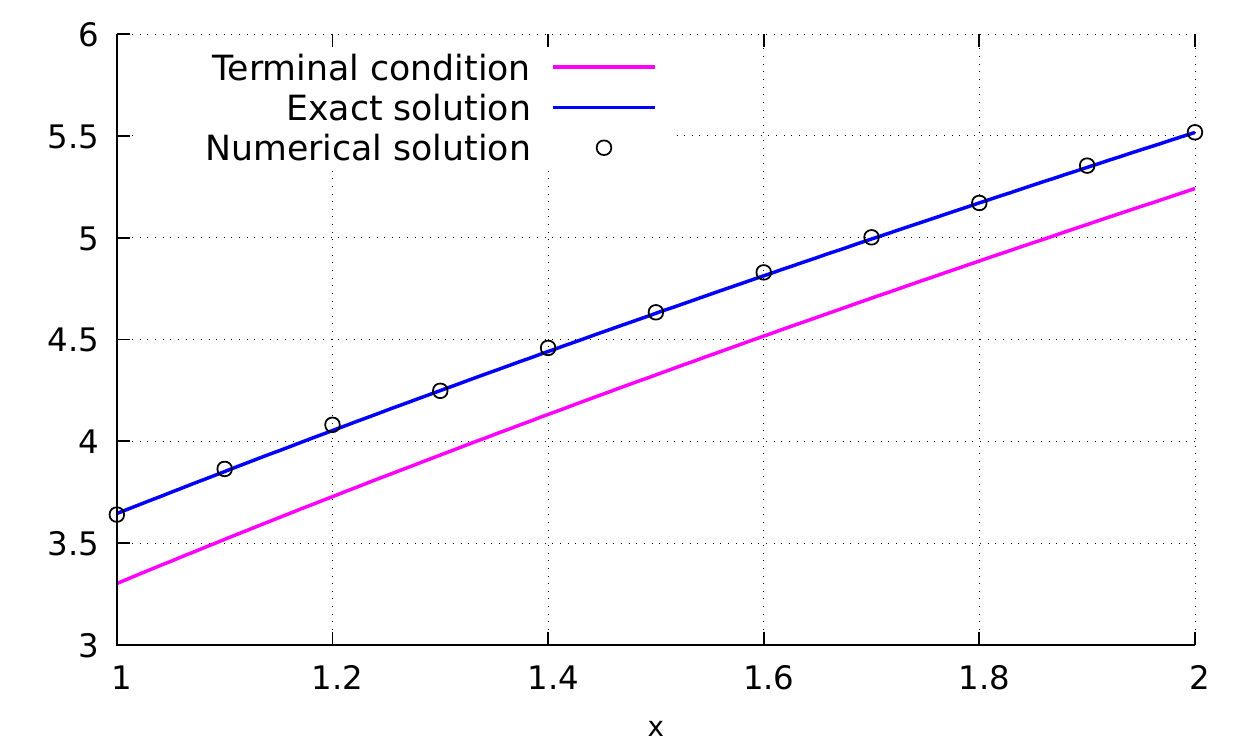}
\vskip-0.1cm
\caption{Dimension $d=1$.}
\end{subfigure}
\begin{subfigure}{.49\textwidth}
\centering
\includegraphics[width=\textwidth]{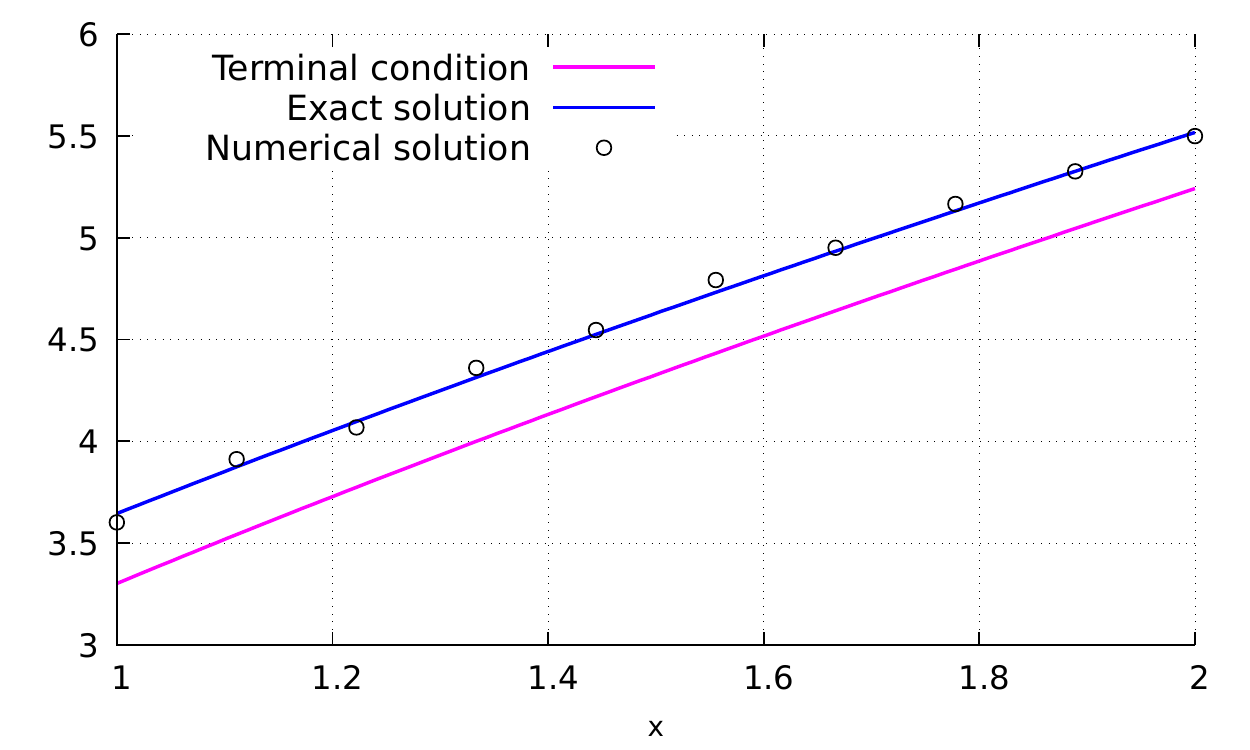}
\vskip-0.1cm
\caption{Dimension $d=5$.}
\end{subfigure}
\caption{Numerical solution $u(0,x)$ of \eqref{dym}.}
\label{f3}
\end{figure}

\vskip-0.3cm

\noindent
As this example and the next one involve derivatives of order greater than one,
they may not be treated by the approach of \cite{labordere} 
due Malliavin weight integrability issues, see page~199 therein. 

\noindent
\item{Example~2-$b)$.}
  For a quasilinear example using non-polynomial nonlinearities,
  consider the equation
\begin{equation}
\label{nonl01} 
\partial_t u(t,x) +
\frac{\alpha}{d} \sum\limits_{i=1}^d \partial_{x_i} u(t,x)
 + \frac{\Delta_x u(t,x)}{1+u^2(t,x)} - 2 u(t,x) = 0, 
\end{equation} 
which admits the traveling wave solution
$$
u(t,x) = \tan \left( \alpha (T-t) + \sum_{i=1}^d x_i \right), 
\qquad (t,x) \in [0,T]\times \real^d, 
$$ 
for $\alpha \in \real$.
In Figure~\ref{f01} we take $T=0.01$, $\alpha =10$ and
one million Monte Carlo samples.
 
\begin{figure}[H]
\centering
\hskip0.2cm
\begin{subfigure}{.49\textwidth}
\hskip0.3cm
\includegraphics[width=\textwidth]{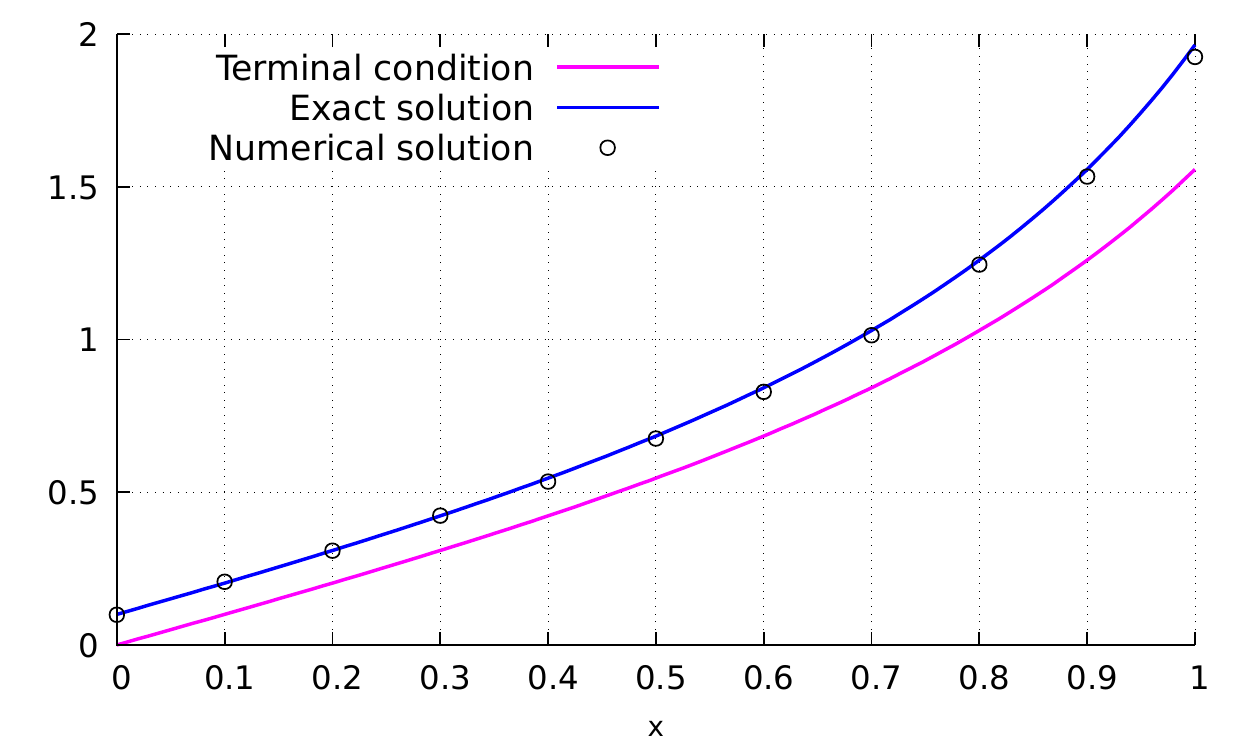}
\vskip-0.1cm
\caption{Dimension $d=1$.}
\end{subfigure}
\begin{subfigure}{.49\textwidth}
\centering
\includegraphics[width=\textwidth]{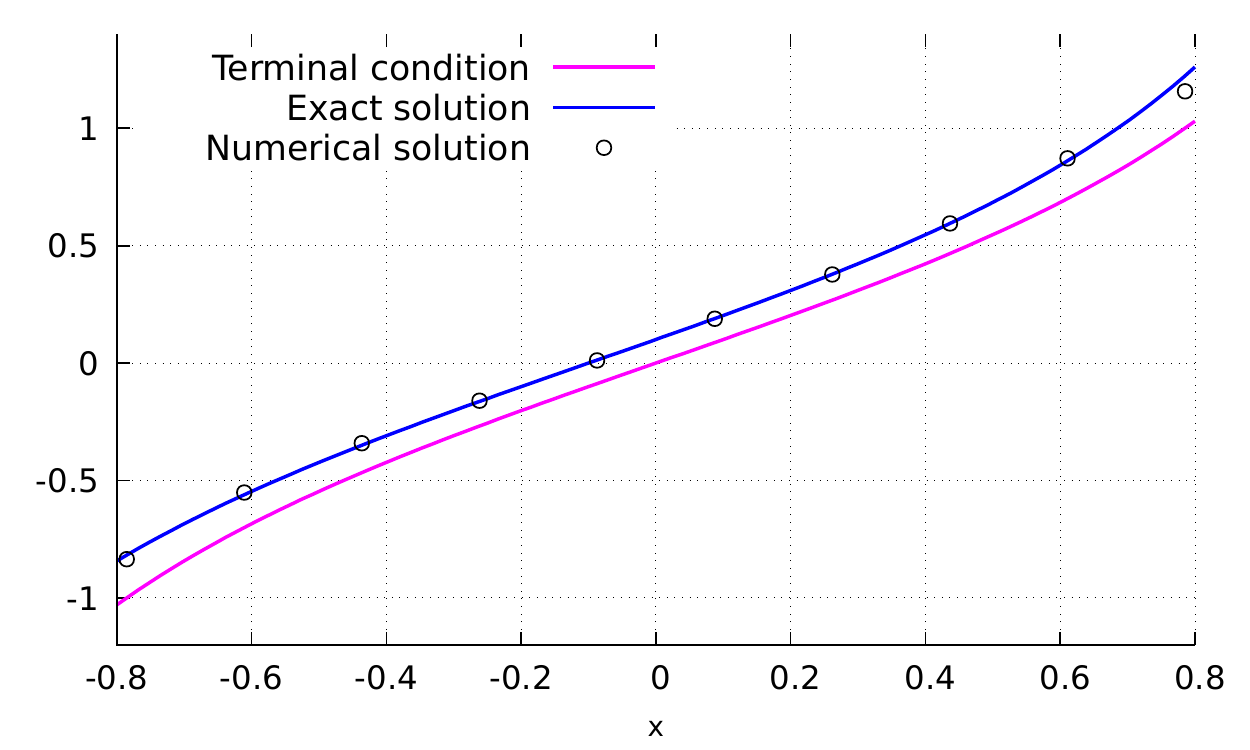}
\vskip-0.1cm
\caption{Dimension $d=5$.}
\end{subfigure}
\caption{Numerical solution $u(0,x)$ of \eqref{nonl01}.}
\label{f01}
\end{figure}

\vskip-0.3cm

\end{description}

\subsubsection*{Quadratic gradient nonlinearity} 
 Consider the Hamilton-Jacobi-Bellman (HJB) equation
\begin{equation}
\label{hjbgl0} 
\partial_t u(t,x) + \Delta_x u(t,x) =
\sum_{i=1}^d ( \partial_{x_i} u(t,x) ) ^2, 
\end{equation} 
 as in \cite{han2018solving} or \S~4.3 of \cite{han2018solvingarxiv},
 with terminal condition $\phi (x) = \log ( (1 + \Vert x\Vert^2 ) /2 )$,
 $x \in \real^d$. 
In Table~\ref{fjl33} we compare our results to the ones obtained
for the estimation of $u(0,x)$ at $x=(0,\ldots , 0)$ for $T=1$ 
with 2000 iterations in Table~2 of \cite{han2018solvingarxiv},
in dimension $d=100$. 
 
\begin{table}[H]
    \centering
\footnotesize 
      \begin{tabular}{|l|c|c|}
   \hline
   & \cite{han2018solvingarxiv} & Coding trees 
   \\
   \hline
 Mean & 4.5977 & 	4.580340
 \\
 \hline
 Standard deviation & 	0.0019 & 	0.001869
  \\
  \hline
  Mean of rel. $L^1$ error & 	0.0017	 & 0.002126
  \\
  \hline
  SD of rel. $L^1$ error & 	0.0004 & 	0.000407
 \\
 \hline
\end{tabular}
\caption{Comparison results with the BSDE method \cite{han2018solvingarxiv} for the HJB equation \eqref{hjbgl0}.} 
\label{fjl33}
\end{table}

\subsubsection*{Fully nonlinear examples} 
In this section, we consider fully nonlinear examples involving 
higher order gradient terms.

\vspace{-0.3cm}

\begin{description}
\item{Example~3-$a)$.}
 For a fully nonlinear example involving a fourth derivative, consider the equation 
 \begin{equation}
   \label{nonl8} 
   \partial_t u(t,x) + \frac{\alpha}{d} \sum\limits_{i=1}^d \partial_{x_i} u(t,x)
   + u(t,x) - \left(\frac{\Delta_x u(t,x)}{12d} \right)^2
   + \frac{1}{d} \sum\limits_{i=1}^d
        \cos \left( \frac{\pi \partial^4_{x_i} u(t,x)}{4!} \right)
   = 0,
\end{equation}
 with terminal condition $\phi (x) := x^4 + x^3 + bx^2 + cx + d$
 where $b = -36/47$, $c = 24b$, $d = 4b^2$, $\alpha = 10$, 
 and solution
 $$
 u(t,x)=\varphi \left( \alpha ( T-t) + \sum\limits_{i=1}^d x_i \right),
 \qquad (t,x) \in [0,T]\times \real^d. 
 $$
 In Figure~\ref{f8} we take $T=0.04$ and $100,000$ Monte Carlo samples. 
 In dimension $d=1$, the graph below is obtained by letting 
 f[y\_\_] :=-y[[3]]/2 + $\alpha$ y[[2]] + y[[1]] - y[[3]]$^2$/144 + Cos[Pi*y[[5]]/24];
 phi[x\_] := x$^4$ + x$^3$ + bx$^2$ + cx + d, 
 and by running Sol[f, 0, T, x, phi, 100000, 4]
 in Mathematica for $x\in [-5, 5]$.

\begin{figure}[H]
\centering
\hskip0.2cm
\begin{subfigure}{.49\textwidth}
\hskip0.3cm
\includegraphics[width=\textwidth]{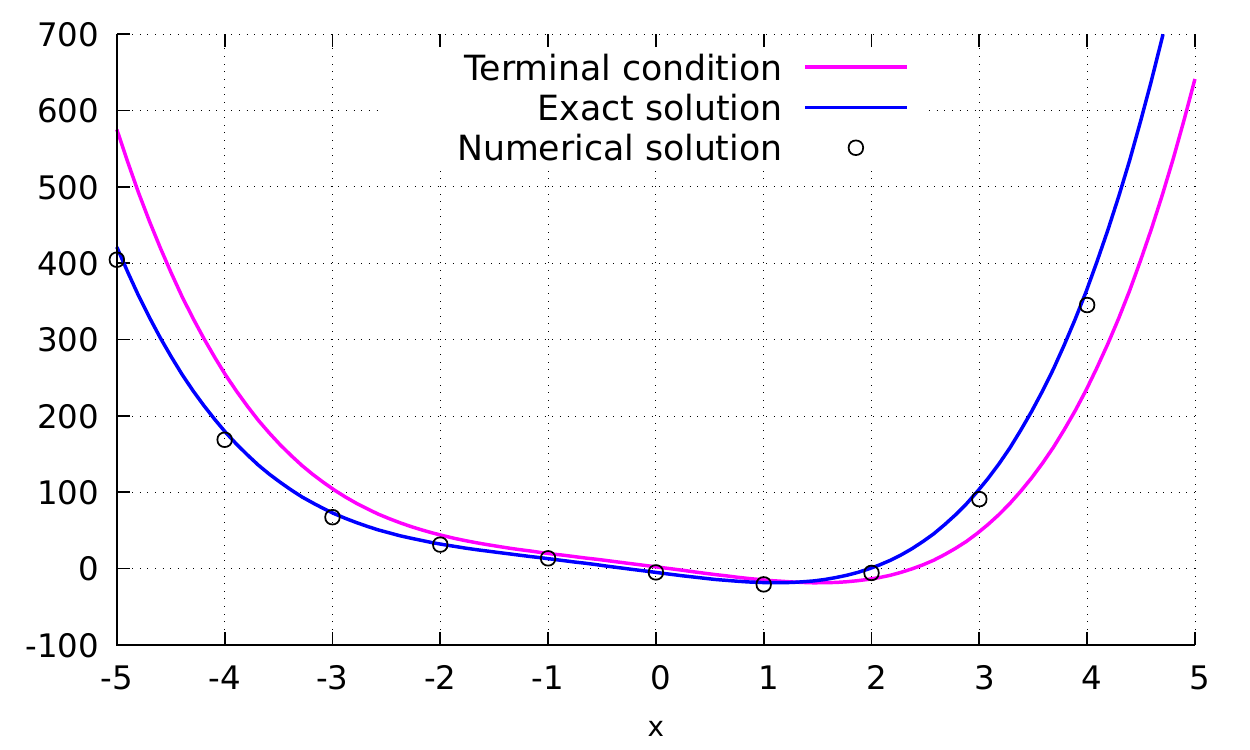}
\vskip-0.1cm
\caption{Dimension $d=1$.}
\end{subfigure}
\begin{subfigure}{.49\textwidth}
\centering
\includegraphics[width=\textwidth]{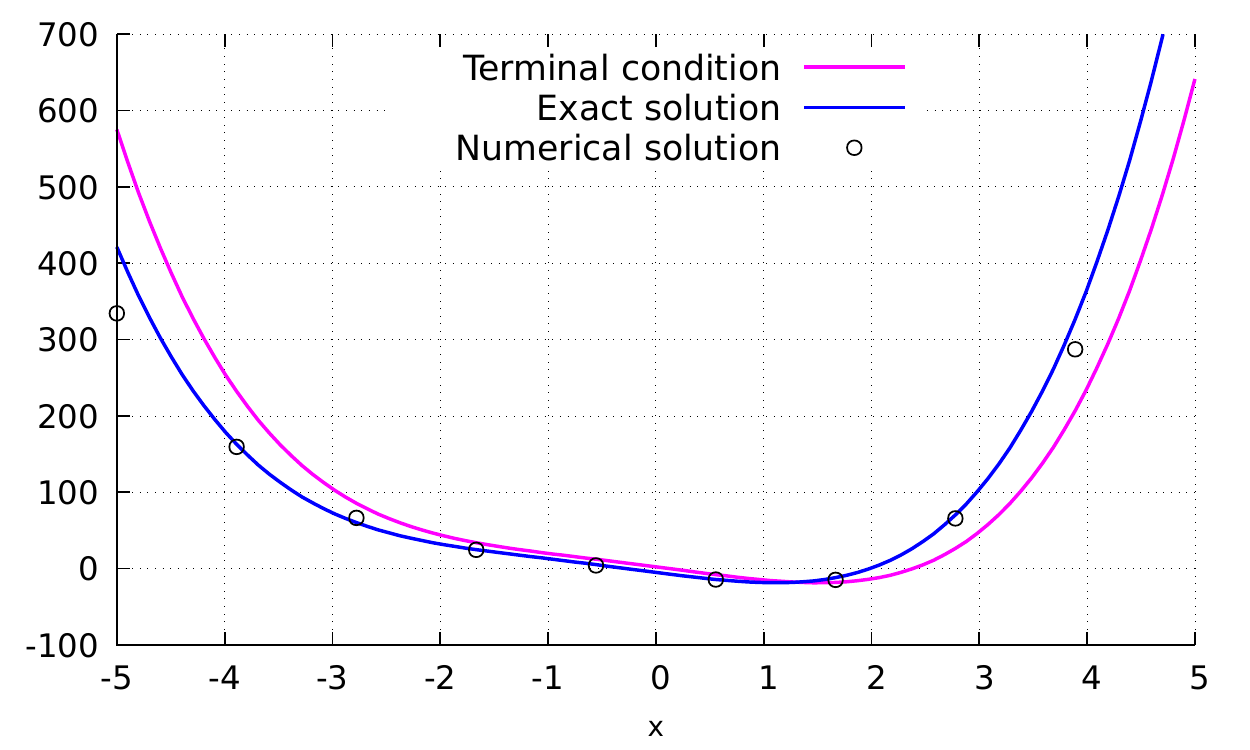}
\vskip-0.1cm
\caption{Dimension $d=5$.}
\end{subfigure}
\caption{Numerical solution $u(0,x)$ of \eqref{nonl8}.}
\label{f8}
\end{figure}

\vspace{-0.6cm}

\item{Example~3-$b)$.}
 For another fully nonlinear example, consider the equation 
 \begin{equation}
   \label{nonl6} 
   \partial_t u(t,x) + \frac{\alpha}{d} \sum\limits_{i=1}^d \partial_{x_i} u(t,x)
   + \log \left(
       \frac{1}{d} \sum\limits_{i=1}^d
       \big(\partial^2_{x_i} u(t,x)\big)^2 + \big(\partial^3_{x_i} u(t,x)\big)^2
   \right)
   = 0,
\end{equation}
 with terminal condition 
$
 \phi (x) = \cos \left( \sum\limits_{i=1}^d x_i \right)$,
 $x \in \real^d$, 
 and solution
$$
u(t,x) = \cos \left( \alpha (T-t) + \sum\limits_{i=1}^d x_i \right),
\qquad (t,x) \in [0,T]\times \real^d, 
$$
 where $\alpha = 5$. 
 In Figure~\ref{f9} we take $T=0.02$ and $100,000$ Monte Carlo samples. 
 In dimension $d=1$, the graph below is obtained by letting 
 f[y\_\_] := $\alpha$*y[[2]] - y[[3]]/2 + Log[y[[3]]$^2$ + y[[4]]$^2$]; 
 phi[x\_] := Cos[x], 
 and by running Sol[f, 0, T, x, phi, 100000, 3]
 in Mathematica for $x\in [-\pi, \pi]$.
 
\begin{figure}[H]
\centering
\hskip0.2cm
\begin{subfigure}{.49\textwidth}
\hskip0.3cm
\includegraphics[width=\textwidth]{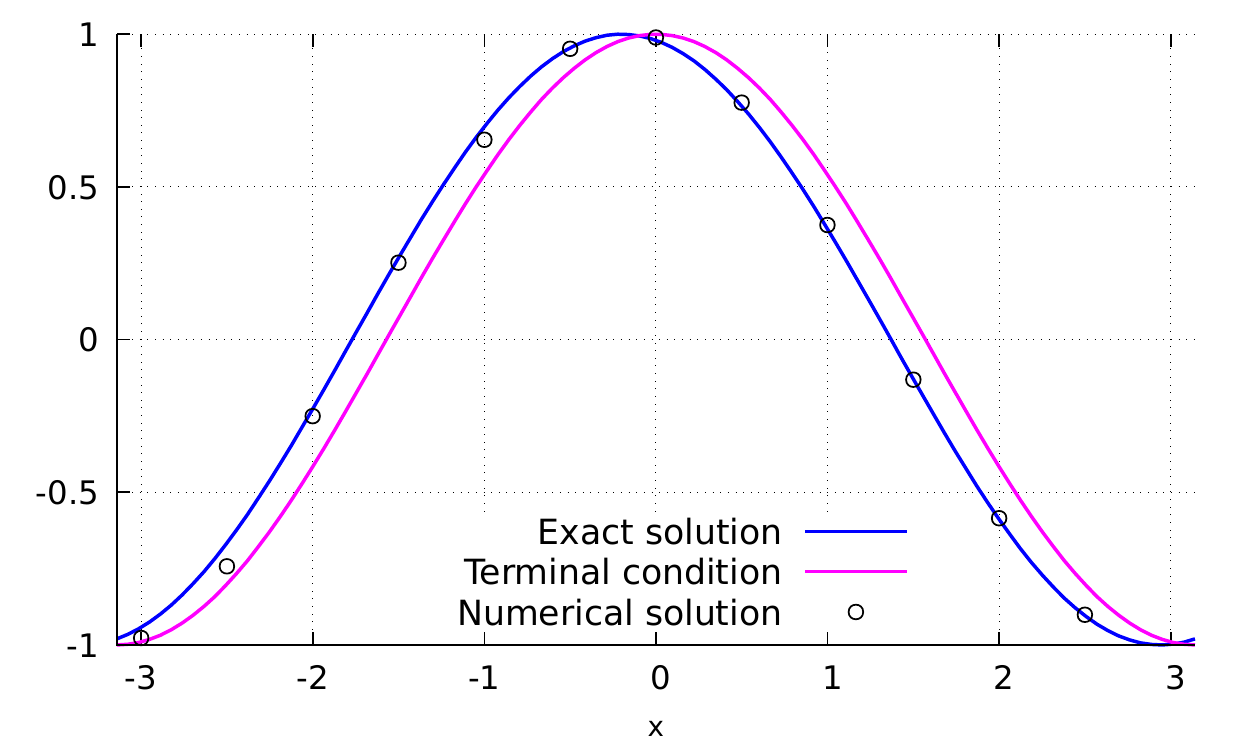}
\vskip-0.1cm
\caption{Dimension $d=1$.}
\end{subfigure}
\begin{subfigure}{.49\textwidth}
\centering
\includegraphics[width=\textwidth]{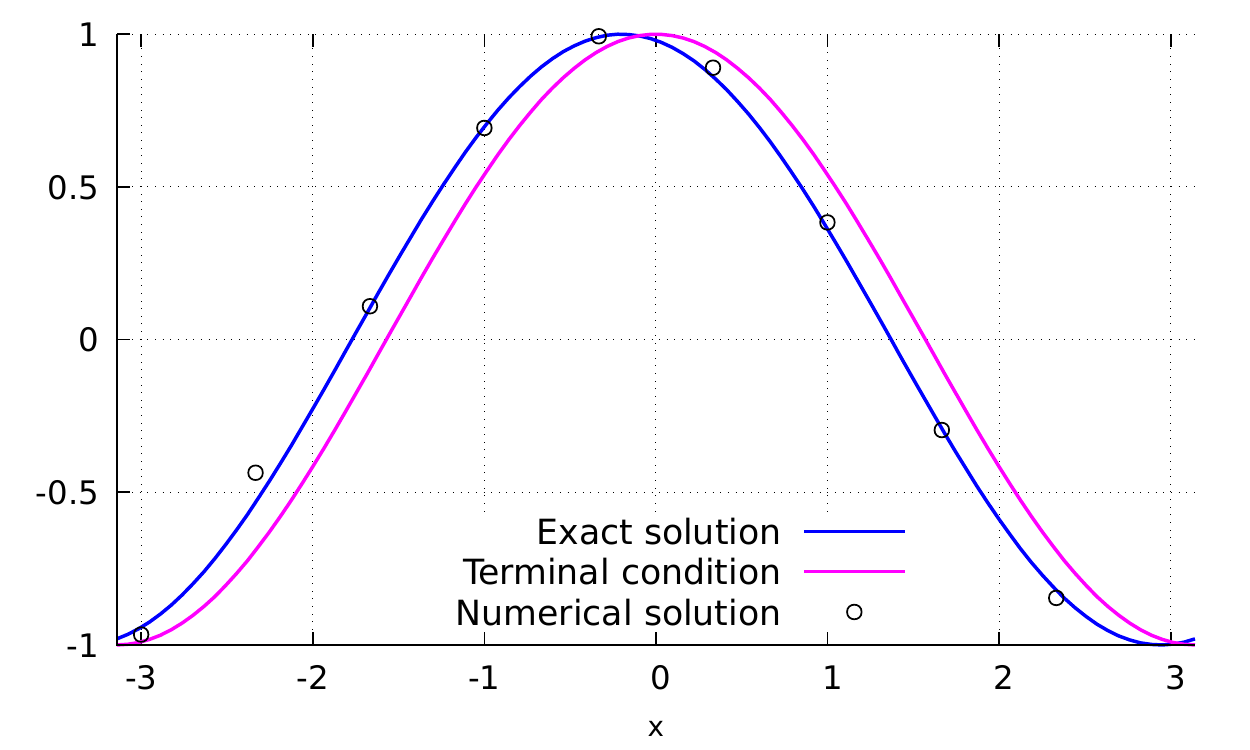}
\vskip-0.1cm
\caption{Dimension $d=5$.}
\end{subfigure}
\caption{Numerical solution $u(0,x)$ of \eqref{nonl6}.}
\label{f9}
\end{figure}

\vspace{-0.6cm}

\end{description} 

\appendix

\section{Computer codes}

The following codes implement the algorithm of Theorem~\ref{t1}
in Mathematica using an exponential distribution $\rho (t) = \re^{-t}$, $t\geq 0$. 
In the above examples the values of fdb[n,k] have been precomputed
by {\em memoization} for $k$ up to $7$ in order to speed up the solution algorithm, where $n$ denotes the highest order of derivative $\partial_x^n$
in \eqref{eq:1}.
The next code implements the mechanism $c \mapsto \mathcal{M}(c)$ 
in the procedure ``codetofunction'' via the combinatorics of the Fa\`a di Bruno
formula written the function ``fdb''.
 
\medskip
 
\begin{lstlisting}[language=Mathematica,caption={2}, title={Fa\`a di Bruno combinatorics and mechanism. 
      } 
  ]
Needs["Combinatorica`"]
kuple[lambda__, s_, l__, k_] := (Module[{m, D, E}, m = Length[lambda]; If[m == 1, Return[Compositions[lambda[[1]], s]]]; E = {}; Do[Do[If[m == 2, D = Append[{k1}, k2], D = Append[k1, k2]]; E = Append[E, D], {k1, kuple[Drop[lambda, -1], s, l, k]}], {k2, Compositions[lambda[[m]], s]}]; Return[E];])
fdb[n_, k_] := fdb[n, k] = (L = {}; Do[Do[Do[ Do[Do[If[n == 1, kv = {ku}, kv = ku]; 
        kw = Sum[kv[[q]], {q, 1, n}]; If[ Signature[l] != 0 && ! MemberQ[kw, 0] && Total[ kw*l] == k , L = Append[L, {k!/Product[l[[p]]!^kw[[p]]*Product[kv[[q]][[p]]!, {q, 1, n}], {p, 1, s}], lambda, kv, l, s}]], {ku, kuple[lambda, s, l, k]}], {l, Select[Tuples[Range[k], s], OrderedQ ]}], {s, 1, k}], {lambda, Compositions[r, n]}], {r, 1, k}]; Return[L])
codetofunction[f_, c__, y0_, phi_] := (Module[{x, y, s1}, If [c == {Id}, Return [phi[y0]]]; If [Length[c] == 1 && c[[1]] < 0, Return [D[phi[y], {y, -c[[1]]}] /. {y -> y0}]]; y = Array[s1, Length[c]]; z = D[f[y], Sequence @@ Transpose[{y, c}]]; Do[z = z /. {s1[k + 1] -> D[phi[x], {x, k}] /. {x -> y0}}, {k, 0, Length[y] - 1}]; Return[z]])
\end{lstlisting}

\noindent
Numerical solution estimates are then computed using the following program,
in which the code $\partial_x^k$ is represented by $\{-k\}$, $k\geq 1$, and
 the code $\big( \partial_{z_0}^{\lambda_0} \cdots \partial_{z_n}^{\lambda_n} f\big)^*
              $ is represented by 
 $\{ \lambda_0 ,  \ldots ,  \lambda_n \} \in \inte^{n+1}$. 
                 
 \medskip

\begin{lstlisting}[language=Mathematica,caption={2}, title={PDE Solution code. 
} 
  ]
  G[x_, tau_] := RandomVariate[NormalDistribution[x, Sqrt[tau]]];
  MCS[f_, t_, tf_, x_, c__, h_, phi_, n_] := (Module[{A, B, U, tau, L, j, l, l1, k1, g, ct}, 
   If[t == tf, Return[phi[x]]]; tau = RandomVariate[ExponentialDistribution[1]];
   If[t + tau >= tf, Return [h* codetofunction[f, c, G[x, tf - t], phi]/Exp[-(tf - t)]]];GS=G[x,tau]; If[c == {Id}, Return[MCS[f, t + tau, tf, GS, ConstantArray[0, n + 1], h/Exp[-tau], phi, n]]]; 
   U = RandomVariate[UniformDistribution[1]][[1]]; 
   If[Length[c] == 1, L = fdb[n + 1, -c[[1]]]; g = Ceiling[U*Length[L]]; A = L[[g]][[1]]* MCS[f, t + tau, tf, GS, L[[g]][[2]], Length[L]*h/Exp[-tau], phi, n]; 
    Do[Do[Do[ A = A*MCS[f, t + tau, tf, GS, {-q - L[[g]][[4]][[j]]}, 1, phi, n], {i1, 1, L[[g]][[3]][[q]][[j]]}], {j, 1, L[[g]][[5]]}], {q, 1, n}]; Return[A]]; 
   l1 = 1 + Sum[Length[fdb[n + 1, k2]], {k2, 1, n}] + (n + 1)^2; 
   If[U <= 1/l1, A = MCS[f, t + tau, tf, GS, ConstantArray[0, n + 1], 1, phi, n]; 
    Return[MCS[f, t + tau, tf, GS, c + UnitVector[n + 1, 1], l1*A*h/Exp[-tau], phi, n]]]; If[U <= (1 + (n + 1)^2)/l1, j = Floor[U*l1*n/(1 + (n + 1)^2)]; 
    l = Ceiling[n*(U - j*(1 + (n + 1)^2)/l1)]; A = MCS[f, t + tau, tf, GS, {-j - 1}, 1, phi, n]; 
    B = MCS[f, t + tau, tf, GS, {-l - 1}, 1, phi, n]; 
    Return[MCS[f, t + tau, tf, GS, c + UnitVector[n + 1, l + 1] + UnitVector[n + 1, j + 1], -l1*A*B*h/Exp[-tau]/2, phi, n]]]; g = Ceiling[U*l1] - 1 - (n + 1)^2; 
   k1 = 1; While[k1 <= n, L = fdb[n + 1, k1]; If[g <= Length[L], Break[]]; g = g - Length[L]; k1++]; 
   A = L[[g]][[1]]* MCS[f, t + tau, tf, GS, L[[g]][[2]], 1, phi, n]; 
   Do[Do[Do[ A = A*MCS[f, t + tau, tf, GS, {-q - L[[g]][[4]][[j]]}, 1, phi, n], {i1, 1, L[[g]][[3]][[q]][[j]]}], {j, 1, L[[g]][[5]]}], {q, 1, n}]; 
   Return[MCS[f, t + tau, tf, GS, c + UnitVector[n + 1, k1 + 1], A*l1*h/Exp[-tau], phi, n]]])
Sol[f_, t_, tf_, x_, phi_, n2_, n_] := (temp = 0; For[i = 1, i <= n2, i++, 
   If[Mod[i, 100000] == 0, Print[i, " sol=", temp/i]]; temp += MCS[f, t, tf, x, {Id}, 1, phi, n];]; Return[temp/n2])
\end{lstlisting}

\footnotesize

\newcommand{\etalchar}[1]{$^{#1}$}
\def\cprime{$'$} \def\polhk#1{\setbox0=\hbox{#1}{\ooalign{\hidewidth
  \lower1.5ex\hbox{`}\hidewidth\crcr\unhbox0}}}
  \def\polhk#1{\setbox0=\hbox{#1}{\ooalign{\hidewidth
  \lower1.5ex\hbox{`}\hidewidth\crcr\unhbox0}}} \def\cprime{$'$}

\end{document}